\UseRawInputEncoding

\documentclass[final]{siamltex}

\usepackage{epsfig}
\usepackage{amsfonts}
\usepackage{amsmath}
\usepackage{amssymb}
\usepackage{verbatim}
\usepackage{sgame}
\usepackage{color}
\usepackage{latexsym}

  \newtheorem{assumption}{\bf{Assumption}}[section]

\begin{document}

\title{Subjective Equilibria under Beliefs of Exogenous Uncertainty for Dynamic Games}

\author{G\"urdal Arslan and Serdar Y\"uksel \thanks{G. Arslan is with the Department of Electrical Engineering,
University of Hawaii at Manoa, 440 Holmes Hall, 2540 Dole Street,
Honolulu, HI 96822, USA. (\email{gurdal@hawaii.edu}). S. Y\"{u}ksel is with the Department of Mathematics and
Statistics, Queen's University, Jeffery Hall, University Avenue,
Kingston, Ontario, CANADA, K7L 3N6. (\email{yuksel@mast.queensu.ca}). This research was
    partially supported by the Natural Sciences and Engineering Research Council of Canada (NSERC).
} }

\maketitle

%\renewcommand{\thefootnote}{\fnsymbol{footnote}}
%
%\footnotetext[2]{G. Arslan is with the Department of Electrical Engineering,
%University of Hawaii at Manoa, 440 Holmes Hall, 2540 Dole Street,
%Honolulu, HI 96822, USA. (\email{gurdal@hawaii.edu}).}
%
%\footnotetext[3]{S. Y\"{u}ksel is with the Department of Mathematics and
%Statistics, Queen's University, Jeffery Hall, University Avenue,
%Kingston, Ontario, CANADA, K7L 3N6. (\email{yuksel@mast.queensu.ca}).}
%
%\renewcommand{\thefootnote}{\arabic{footnote}}

%\slugger{mms}{xxxx}{xx}{x}{x--x}%slugger should be set to mms, siap, sicomp, sicon, sidma, sima, simax, sinum, siopt, sisc, or sirev

%%%%%%%%%%%%%%%%%%%%%%%%%%%%%%%%%%%%%%%%%%%%%%%%%%%%%%%%%%%%%%%%%%%%%%%%%%%%%%%%

\begin{abstract}
We present a subjective equilibrium notion (called "subjective equilibrium under beliefs of exogenous uncertainty (SEBEU)" for stochastic dynamic games in which each player chooses her decisions under the (incorrect) belief that a stochastic environment process driving the system is exogenous whereas in actuality this process is a solution of closed-loop dynamics affected by each individual player. Players observe past realizations of the environment variables and their local information. At equilibrium, if players are given the full distribution of the stochastic environment process as if it were an exogenous process, they would have no incentive to unilaterally deviate from their strategies. This notion thus generalizes what is known as the static price-taking behavior in prior literature to a stochastic and dynamic setup. We establish existence of SEBEU, study various properties and present explicit solutions. We obtain the $\epsilon$-Nash equilibrium property of SEBEU when there are many players.
\end{abstract}

%\begin{keywords}
%Game theory, subjective equilibrium, distributional consistency, stochastic dynamic games.
%\end{keywords}

%\begin{AMS}
%93E03, 93E35, 91A26, 68T05
%\end{AMS}

%\pagestyle{myheadings}
%\thispagestyle{plain}
%\markboth{SUBJECTIVE EQUILIBRIA UNDER BELIEFS OF EXOGENOUS UNCERTAINTY}{}

%%%%%%%%%%%%%%%%%%%%%%%%%%%%%%%%%%%%%%%%%%%%%%%%%%%%%%%%%%%%%%%%%%%%%%%%%%%%%%%%
\section{Introduction}
In this paper, we consider a non-cooperative stochastic multi-stage game with finite and infinite number of players and discrete-time state dynamics. Each player's decision at any stage is used to control the stochastic evolution of the player's own state dynamics. Each player's objective is to minimize her own expected total cost, which sums up the stage-cost incurred by the player at each stage.
The stage cost and the state evolution of each player is influenced by the so-called environment variables in addition to the player's own decisions and states. The environment variables are the same across all players and are determined as a function of all player's decisions.
Therefore, players' problems are coupled only through the environment variables, i.e., they are otherwise decoupled.

Players in our model know their own individual parameters such as their set of strategies, cost functions, discount factors, and have access to the history of their own state and control variables as well as the environment variables. However, players neither have any information about the other players nor do they know how the environment variables are generated.
Players make the modeling assumption that the environment variables are generated by an independent exogenous random process with a known probability distribution. Based on such incorrect modeling assumptions, players choose their strategies to minimize their own long-term cost. Our main concern in this paper is to find player strategies which (i) are optimal with respect to the players' assumed models for the environment variables, and (ii) generate environment variables whose probability distribution is consistent with the probability distributions assumed by all players for the environment variables.
Such strategies constitute a form of subjective equilibrium at which players cannot detect the falseness of their modeling assumptions unless they actively probe: the environment variables generated at a subjective equilibrium are distributionally consistent with their beliefs when viewed in isolation as far as the marginal probability measure is concerned, hence they would not have an incentive to deviate to an alternative strategy. At such an equilibrium, to be called a \textit{subjective equilibrium under beliefs of exogenous uncertainty (SEBEU)}, if players are given the full distribution of the stochastic environment process as if it is an independent and exogenous process (that is, primitive, generated by nature independently), they would have no incentive to unilaterally deviate from their equilibrium strategies.

We must emphasize that, in games with finite set of players, the notion of subjective equilibrium under beliefs of exogenous uncertainty is different from that of Nash equilibrium and its various refinements. At a Nash equilibrium, each player's beliefs of the environment is objective (in particular, a player does not view the environment variables as exogenous), hence, each player takes into account not only the direct effect of her decisions on her cost but also the indirect (i.e., the closed-loop) effect of her decisions on her cost through the environment variables. This indirect effect through the environment variables is not taken into account at an SEBEU at which the environment variables are viewed as exogenous.
%\sy{Even though, in a Nash equilibrium, the environment variables may be induced by the policies of the decision makers given, each individual decision maker's policy affect the induced cost both through the immediate actions/costs, but also the environment variables in the closed-loop through the policies of the other players. As a particularly insightful case in point; consider an identical interest (team) setup involving Linear Quadratic Gaussian setups. In this case, a price taking equilibrium will be linear, however, the optimal decentralized stochastic control policies will not be linear as the optimization problem will no longer be convex in the individual decision maker's policies.}

The related notion of \textit{price-taking behavior} is well known in mathematical economics. For example, the fundamental notion of competitive equilibrium at which the players in an exchange economy optimize their utilities with respect to market clearing prices hinges on the assumption that the agents take the prices as given \cite{arrow1954existence}, \cite{radner1968competitive}.
Although the notion of price taking behavior has been used in different settings in the literature (e.g., the notion of equilibrium used in \cite{kelly1997charging,kelly1998rate} for communication networks), to the best of authors' knowledge, the notion of subjective equilibrium introduced in this paper in its full generality (with distributional consistency) does not exist in the current literature.

\cite{green1980noncooperative} introduced a notion of subjective equilibrium, similar to the notion of subjective equilibrium introduced in this paper, for dynamic markets with no state dynamics and a restrictive set of strategies. The results of \cite{green1980noncooperative} address the questions of (i) whether Nash equilibria include price-taking equilibria when the set of firms is a continuum, and (ii) whether Nash equilibria in a sequence of increasingly large finite economies approach price-taking equilibria of a non-atomic limiting economy, in particular, they do not address the existence of subjective equilibrium in economies with a finite set of firms. Similar issues of approximating Nash equilibria in large finite games by the equilibria of non-atomic limiting games have been dealt with in a number of references under the general theme of large anonymous games; see for example \cite{jovanovic1988anonymous} and \cite{CARMONA2020105015} and the references therein.

The notion of mean-field equilibrium \cite{huang2006,lasry2007mean} has also a similar flavor in that the mean-field equilibria are considered as the limit of the Nash equilibria as the number of players tend to infinity.
Sufficient conditions for the convergence of Nash equilibria to mean-field equilibria (possibly along subsequences), as the number of (nearly) identical players tend to infinity, have been presented in \cite{fischer2017connection, lacker2018convergence, bardi2014linear, lasry2007mean, arapostathis2017solutions,sanjari2018optimal,sanjari2019optimal,sanjari2020optimality}, with the last three papers focusing on identical interest games or stochastic teams. We refer interested readers to \cite{carmona2018probabilistic,Caines2017} for a literature review and a detailed summary of some recent results on mean-field games.

In contrast to these references where the prevailing concept is that of Nash equilibrium, we concern ourselves in this paper with a notion of subjective equilibrium (which generally does not coincide with Nash equilibrium when the set of players is finite) in which the players ignore the effects of their strategies on the environment variables (e.g., prices) due to their beliefs of exogenous uncertainty in the environment variables.

Other relevant references include \cite{adlakha2015equilibria, weintraub2011industry,wiszniewska2017redefinition,dudebout2012empirical,dudebout2014exogenous,Kalai:GEB:1995,wiszniewska2014open}. The notion of equilibrium introduced in \cite{adlakha2015equilibria}, called stationary equilibrium, is based on the assumption that the players view the population states (which play the role of the environment variables in this paper) as deterministic and constant. This equilibrium notion requires that the constant population states assumed by the players equal the long-run average of the actual population state (averaged over infinite number of players). If the number of players is finite or the population state is not initially at steady-state, a player observing the actual population state can clearly reject the hypothesis of a constant deterministic population state.

The notion of oblivious equilibrium utilized in \cite{weintraub2011industry} is completely analogous to the stationary equilibrium in \cite{adlakha2015equilibria}.
The notion of belief distorted Nash equilibrium introduced in \cite{wiszniewska2017redefinition} for a deterministic system is closer in spirit to the subjective equilibrium in this paper and it requires players' (probabilistic) beliefs to concentrate on the actual deterministic play path. The notion of empirical evidence equilibrium introduced in \cite{dudebout2012empirical,dudebout2014exogenous} for stochastic games (with finite state and control spaces) is also conceptually similar to the equilibrium notion in this paper and it requires an eventual consistency between the objective and the subjective probabilities of the signal variable (environment variable in this paper) conditioned on the assumed Markovian state for the signals. There is also a separate body of work in the literature on dynamic games with a continuum of players, which includes private as well as global state variables (as in our paper); see \cite{wiszniewska2014open} and the references therein.

We also cite \cite{Kalai:GEB:1995} as a related key reference. The subjective equilibrium introduced in \cite{Kalai:GEB:1995} and the one in this paper are distinct concepts due to the differences in the beliefs held by the players with regard to the nature of the environment variables and in the requirements of consistency of the beliefs with the actual play.
While the consistency of the beliefs on the play paths, which include correlations between player strategies and the environment variables, are required
in \cite{Kalai:GEB:1995}, players in this paper ignore such correlations, which is arguably more appropriate in certain applications, for example, the Autonomous Demand Response application on price-sensitive control of energy usage in multiple (possibly very large number of) units \cite{constantopoulos1991estia} \cite{roozbehani2012volatility}; see Subsection~\ref{ss:ex} for a simplified single-stage example of such an application.

{\bf Contributions.} (i) We introduce a subjective equilibrium concept with distributional consistency, in the sense of probability measures on sequence space, for stochastic dynamic games. As reviewed and motivated in detail above, while related concepts of equilibria have appeared in the literature, prior studies have only focused on the more restrictive static or deterministic setups. In particular, our main contribution is to introduce the concept of {\it subjective equilibrium under beliefs of exogenous uncertainty (SEBEU)}, which is applicable to a very general class of problems. (ii) We obtain a general existence result on SEBEU in mixed strategies for compact metric space models with finite as well as infinite time horizons, which include all finite space models (where the state, control, and environment variables belong to finite sets) as a special case. (iii) We establish the near person-by-person optimality of an SEBEU (which constitutes an approximate Nash equilibrium) when there are large number of ``small players''  and the environment variables are generated by certain aggregate statistics on player states-decisions. This result on near person-by-person optimality is also extended to pure strategies under additional assumptions which, in particular, requires identical players. (iv) When the number of decision makers is infinite, we establish that an SEBEU is also a Nash equilibrium.
We emphasize that the multi-stage nature of the analysis and the modeling of the environment variables as an exogenous stochastic process, with possibly arbitrary correlations between the environment variables, in the closed-loop lead to technical challenges to be addressed in the paper.

{\bf Notation.} $\mathbb{N}_0$ and $\mathbb{N}$ denote the nonnegative and positive integers, respectively;  $\mathbb{R}$ denotes the real numbers; %$[\cdot]^{\prime}$ denotes the transpose of a vector or matrix;
$|\cdot|$ denotes the Euclidean norm for a vector;
%$|x|_Q^2:=x^{\prime}Qx$ for a vector $x$ and a matrix $Q$;
%$A\succeq0$ denotes a nonnegative definite matrix $A$, whereas $A\succ0$ denotes a positive definite matrix $A$; $A\succeq B$ denotes $A-B\succeq0$;
$\mathcal{P}(\mathbb{T})$ denotes the set of probability measures on the Borel sigma algebra $\mathcal{B}(\mathbb{T})$ of a topological space $\mathbb{T}$; $P[\cdot]$, $E[\cdot]$, and $\textrm{cov}[\cdot]$ denote the probability, the expectation, and the covariance respectively ($P^{\mu}[\cdot]$, $E^{\mu}[\cdot]$, and $\textrm{cov}^{\mu}[\cdot]$ are also used to emphasize the underlying strategy $\mu$);
%$X\sim\mathcal{N}(\theta,\Sigma)$ denotes a Gaussian random vector with mean $\theta$ and covariance $\Sigma$;
$m_{\mathbb{M}}$ denotes a metric on a metric space $\mathbb{M}$\footnote{ Moreover, when we speak of a Cartesian product of a countable collection of metric spaces
$((\mathbb{M}_t,m_{\mathbb{M}_t}))_{t\in[0,T)}$, we assume that the product space $\mathbb{M}:=\times_{t\in[0,T)} \mathbb{M}_t$ is equipped with the metric
\begin{equation}
\label{eq:pmet}
m_{\mathbb{M}}(a,b) := \sum_{t\in[0,T)} 2^{-t} \frac{m_{\mathbb{M}_t}(a_t,b_t)}{1+m_{\mathbb{M}_t}(a_t,b_t)}
\end{equation}
for all $a=(a_t)_{t\in[0,T)},b=(b_t)_{t\in[0,T)}\in\mathbb{M}$. The metric $m_{\mathbb{M}}$, known as the product metric, generates the product topology \cite{aliprantis2006infinite}.}.

%%%%%%%%%%%%%%%%%%%%%%%%%%%%%%%%%%%%%%%%%%%%%%%%%%%%%%%%%%%%%%%%%%%%%%%%%%%%%%%%%
%\section{MODEL}
%\section{Model}
%\label{se:model}

\section{Subjective equilibrium under beliefs of exogenous uncertainty}
\subsection{Model and definitions}\label{se:model}
Consider a decentralized stochastic system with $N$ decision makers, $N\in\mathbb{N}$, where the $i$-th decision maker is referred to as DM$^i$, $i\in [1,N]:=\{1,\dots,N\}$. Time variable $t\in[0,T):=\{0,\dots,T-1\}$ is an integer where $T\in\mathbb{N}$ is the time horizon; we will also allow for $T= \infty$ in which case $[0,T) = \mathbb{N}_0$. Each DM$^i$ has its own state $x_t^i\in\mathbb{X}^i$, control input $u_t^i\in\mathbb{U}^i$, and random disturbance $w_t^i\in\mathbb{W}^i$ at time $t\in[0,T)$. In addition, there are so-called environment variables denoted by $y_t\in\mathbb{Y}$ at time $t\in[0,T)$. Each DM$^i$'s state evolves as
\begin{equation}
\label{eq:x}
x_{t+1}^i =f_t^i(x_t^i,u_t^i,y_t,w_t^i), \qquad t\in [0,T)
\end{equation}
starting from some (possibly random) initial state $x_0^i\in\mathbb{X}^i$, where $f_t^i:\mathbb{X}^i\times\mathbb{U}^i\times\mathbb{Y}\times\mathbb{W}^i \to \mathbb{X}^i$ is some function.
The environment variable $y_t$ is produced by
\begin{align}
\label{eq:p}
y_t =  & g_t(x_t^0, x_t, u_t,\xi_t), \qquad t\in[0,T) \\
\label{eq:y}
x_{t+1}^0 = & f_t^0(x_t^0, x_t, u_t, w_t^0), \qquad t\in[0,T)
\end{align}
where $x_t^0\in\mathbb{X}^0$ is the state of the environment dynamics at time $t\in[0,T)$ ($x_0^0\in\mathbb{X}^0$ is some possibly random initial state),
$x_t=(x_t^1,\dots,x_t^N)\in\mathbb{X}:=\mathbb{X}^1\times\cdots\times\mathbb{X}^N$ is the state of all DMs at time $t\in[0,T)$,
$u_t=(u_t^1,\dots,u_t^N)\in\mathbb{U}:=\mathbb{U}^1\times\cdots\times\mathbb{U}^N$ is the joint control input at time $t\in[0,T)$,  $\xi_t\in\Xi$, $w_t^0\in\mathbb{W}^0$  are random disturbances at time $t\in[0,T)$, and $g_t:\mathbb{X}^0\times\mathbb{X}\times\mathbb{U}\times\Xi \to \mathbb{Y}$, $f_t^0:\mathbb{X}^0\times\mathbb{X}\times\mathbb{U}\times\mathbb{W}^0 \to \mathbb{X}^0$ are some functions, $t\in[0,T)$.
We assume throughout that each of the following spaces is nonempty and metric:
\begin{equation}
\label{eq:sets}
\mathbb{X}^0,\dots,\mathbb{X}^N,\mathbb{U}^1,\dots,\mathbb{U}^N,\mathbb{W}^0,\dots,\mathbb{W}^N,\mathbb{Y},\Xi.
\end{equation}

Throughout the paper, each DM$^i$ selects its control input $u_t^i$ at time $t$ based on the history of its observations until time $t$, i.e., each DM$^i$ uses a closed-loop strategy in the language of optimal control theory as opposed to using an open-loop strategy which would be a function of the time variable $t$ alone. More precisely, each DM$^i$ has access to the information $(I_t^i,Y_{t-1})$ at time $t\in[0,T)$, where
$$I_0^i:=x_0^i, \ \ Y_{-1}:=\emptyset, \ \ I_t^i   := ((x_k^i)_{k\in[0,t]},(u_k^i)_{k\in[0,t)}), \ \ Y_{t-1} :=(y_k)_{k\in[0,t)}, \quad t\in[1,T).$$
We note that $I_t^i$ denotes DM$^i$'s private information  whereas $Y_{t-1}$ denotes the common information available to all DMs at time $t$.
A \textit{pure strategy} for DM$^i$ is a collection $s^i=(s_t^i)_{t\in[0,T)}$  of mappings $s_t^i:\mathbb{I}_t^i \times \mathbb{Y}^{[0,t)} \to \mathbb{U}^i$, where $\mathbb{I}_t^i:=(\mathbb{X}^i)^{[0,t]}\times(\mathbb{U}^i)^{[0,t)}$,
$t\in[0,T)$. For both $T\in\mathbb{N}$ and $T=\infty$, $(I_T^i,Y_{T-1})\in\mathbb{I}_T^i\times\mathbb{Y}^{[0,T)}$ denotes DM$^i$'s entire history of observations.
We let $\mathbb{S}^i$ denote the set of pure strategies for each DM$^i$. If DM$^i$ employs a pure strategy $s^i\in\mathbb{S}^i$, then
\begin{equation}
\label{eq:u}
u_t^i=s_t^i(I_t^i,Y_{t-1}), \qquad t\in[0,T).
\end{equation}
For any \textit{joint pure strategy} $s=(s^1,\dots,s^N)\in\mathbb{S}:=\mathbb{S}^1\times\cdots\times\mathbb{S}^N$ employed,
each DM$^i$ incurs a long-term cost
\begin{equation}
\label{eq:ltc1}
\bar{J}^i(s^i,s^{-i}) = E \Bigg[\sum_{t\in[0,T)} (\beta^i)^t c_t^i(x_t^i,u_t^i,y_t) + (\beta^i)^T c_T^i(x_T^i) \Bigg]
\end{equation}
where $s^{-i}\in\mathbb{S}^{-i}:=\times_{j\not=i} \mathbb{S}^j$ denotes the pure strategies of all DMs other than DM$^i$, $\beta^i\in(0,1]$ is DM$^i$'s discount factor ($\beta^i\in(0,1)$ when $T=\infty$),
$c_t^i:\mathbb{X}^i\times\mathbb{U}^i\times\mathbb{Y}\to\mathbb{R}$ is some function which determines DM$^i$'s cost at each time $t\in[0,T)$, and $c_T^i:\mathbb{X}^i\to\mathbb{R}$ determines DM$^i$'s terminal cost ($c_T^i\equiv0$ when $T=\infty$).
The expectation in (\ref{eq:ltc1}) is taken with respect to the probability distribution over the set $\mathbb{I}_T^i\times\mathbb{Y}^{[0,T)}$ of histories  induced by $s\in\mathbb{S}$.
Note that each DM$^i$'s cost is influenced by the other DMs only through the environment variables $Y_{T-1}$.

Suppose now that each DM$^i$ views the sequence of environment variables as an {\bf independent exogenous} (that is, primitive, generated by nature independently of the all other primitive random variables) random sequence $Z_{T-1}:=(z_t)_{t\in[0,T)}$ with a given probability distribution $\zeta\in\mathcal{P}(\mathbb{Y}^{[0,T)})$  instead of the endogenous sequence $Y_{T-1}$ generated by (\ref{eq:x})-(\ref{eq:y}). We note that $Z_{T-1}$ need not be a collection of independent random variables. We define DM$^i$'s cost corresponding to a pure strategy $s^i\in\mathbb{S}^i$ and a deterministic sequence of environment variables  $\bar{Z}_{T-1}=(\bar{z}_t)_{t\in[0,T)}\in \mathbb{Y}^{[0,T)}$ as
\begin{equation}
\label{eq:ltc2}
J^i(s^i,\bar{Z}_{T-1}) := E \Bigg[\sum_{t\in[0,T)} (\beta^i)^t c_t^i\big(x_t^i,u_t^i,\bar{z}_t\big)   + (\beta^i)^T c_T^i(x_T^i) \Bigg]
\end{equation}
where
\begin{align*}
x_{t+1}^i = & f_t^i(x_t^i,u_t^i,\bar{z}_t,w_t^i), \qquad t\in[0,T) \\
u_t^i = & s_t^i(I_t^i,\bar{Z}_{t-1}), \qquad t\in[0,T)
\end{align*}
with $\bar{Z}_{-1}:=\emptyset$, and $\bar{Z}_{t-1} :=(\bar{z}_0,\dots,\bar{z}_{t-1})$, for $t\in[1,T)$. Note that the expectation in (\ref{eq:ltc2}) is taken with respect to the probability distribution over the set $\mathbb{I}_T^i$ of DM$^i$'s state-control histories  induced by $(s^i,\bar{Z}_{T-1})\in\mathbb{S}^i \times \mathbb{Y}^{[0,T)}$.
Accordingly, each DM$^i$ who models the environment variables as an independent exogenous random sequence with the probability distribution $\zeta\in\mathcal{P}(\mathbb{Y}^{[0,T)})$ aims to minimize $E^{\zeta} [J^i(s^i,\cdot)]$ by choosing a pure strategy $s^i\in\mathbb{S}^i$.
Therefore, from DM$^i$'s viewpoint, no DM has any influence on the environment variables and the cost $E^{\zeta} [J^i(s^i,\cdot)]$ is independent of the strategies of all DMs other than DM$^i$.

In all cases considered in the paper, the appropriate conditions will be imposed to ensure that the long-term cost (\ref{eq:ltc1}) and (\ref{eq:ltc2}) are well-defined for any joint strategy $s\in\mathbb{S}$ and independent exogenous distribution $\zeta\in\mathcal{P}(\mathbb{Y}^{[0,T)})$ for the environment variables.
\\

%\begin{definition} Let $s=(s^1,\dots,s^N)\in\mathbb{S}$ be a joint strategy and $\zeta_{s}\in\mathcal{P}(\mathbb{Y}^{[0,T)})$ be the probability distribution of the environment variables generated by $s$ endogenously through (\ref{eq:x})-(\ref{eq:u}). The joint strategy $s$ is called a {\it subjective equilibrium under beliefs of exogenous uncertainty (SEBEU) in pure strategies} if
%\begin{equation}
%\label{eq:br}
%E^{\zeta_s} [J^i(s^i,\cdot)]  \leq E^{\zeta_s} [ J^i(\tilde{s}^i,\cdot) ], \quad \forall  i\in[1,N], \tilde{s}^i\in\mathbb{S}^i
%\end{equation}
%where the environment variables are (incorrectly) assumed, by each DM$^i$, to be an independent exogenous random sequence with the probability distribution $\zeta_{s}$.\\
%\end{definition}

\begin{definition}
\label{def:sebeu}
A joint strategy $s=(s^1,\dots,s^N)\in\mathbb{S}$ is called a {\it subjective equilibrium under beliefs of exogenous uncertainty (SEBEU) in pure strategies} if there exists subjective beliefs $\zeta^1,\dots,\zeta^N\in\mathcal{P}(\mathbb{Y}^{[0,T)})$, where $\zeta^i$ is DM$^i$'s belief about the environment variables, such that the following two conditions hold.
\begin{itemize}
\item[(i)] Each DM$^i$'s strategy $s^i$ is optimal with respect to its subjective belief $\zeta^i$, i.e.,
$$E^{\zeta^i} [J^i(s^i,\cdot)]  \leq E^{\zeta^i} [ J^i(\tilde{s}^i,\cdot) ], \quad \forall  i\in[1,N], \tilde{s}^i\in\mathbb{S}^i$$
\item[(ii)] Each DM$^i$'s subjective belief $\zeta^i$ is consistent with the objective probability distribution $\zeta_{s}$ of the environment variables, i.e.,
$$\zeta^i(B) = \zeta_s(B), \quad \forall  i\in[1,N], B\in\mathcal{B}(\mathbb{Y}^{[0,T)})$$
where $\zeta_s\in\mathcal{P}(\mathbb{Y}^{[0,T)})$ is generated by $s$ endogenously through (\ref{eq:x})-(\ref{eq:u}).\\
\end{itemize}
\end{definition}

Note that, in an SEBEU, each DM$^i$ ignores the influence of its own strategy $s^i$ on the environment variables in minimizing its long-term cost. In contrast, the well-known concept of {\it Nash equilibrium in pure strategies} requires a joint pure strategy $s=(s^1,\dots,s^N)\in\mathbb{S}$ to satisfy
\begin{equation}
\bar{J}^i(s^i,s^{-i}) \leq \bar{J}^i(\tilde{s}^i,s^{-i}), \quad \forall i\in[1,N], \tilde{s}^i\in\mathbb{S}^i
\label{eq:neq}
\end{equation}
where each DM$^i$ takes into account the entire influence of its own strategy $s^i$ on its long-term cost including through the environment variables.

The concepts introduced above are naturally extended to so-called mixed strategies which are probability distributions over the pure strategies. A \textit{mixed strategy} for DM$^i$  is a probability distribution $\gamma^i\in\mathcal{P}(\mathbb{S}^i)$ over DM$^i$'s set of pure strategies.
We let $\Gamma^i:=\mathcal{P}(\mathbb{S}^i)$ denote the set of mixed strategies for each DM$^i$. If DM$^i$ employs a mixed strategy $\gamma^i\in\Gamma^i$, then DM$^i$ makes an independent random draw $s^i\in\mathbb{S}^i$ according to $\gamma^i$ and employs the pure strategy $s^i$. A joint mixed strategy $\gamma=(\gamma^1,\dots,\gamma^N)\in \Gamma:=\Gamma^1\times\cdots\Gamma^N$ induces a probability distribution over the set $\mathbb{Y}^{[0,T)}$ of environment variables defined as: for any $B\in\mathcal{B}(\mathbb{Y}^{[0,T)})$,
\begin{equation}
\label{eq:zeta_gamma}
\zeta_{\gamma}(B):=E^{\gamma} [\zeta_s(B)]
\end{equation}
where $\zeta_s\in\mathcal{P}(\mathbb{Y}^{[0,T)})$ is the probability distribution of the environment variables generated by $s\in\mathbb{S}$ endogenously via (\ref{eq:x})-(\ref{eq:u}).
The long-term cost of each DM$^i$ under pure strategies is extended to mixed strategies as $E^{\gamma} [ \bar{J}^i ]$ or $E^{(\gamma^i,\zeta)} [J^i]$ (for an independent exogenous distribution $\zeta\in\mathcal{P}(\mathbb{Y}^{[0,T)})$ for the environment variables). Using these expected long-term cost for each DM$^i$, the notion of an SEBEU or a Nash equilibrium extends to mixed strategies in the most natural way. Finally, the pure strategies in each $\mathbb{S}^i$ can be identified with the extreme elements of $\Gamma^i$ each of which assigns probability one to a single element of $\mathbb{S}^i$.

%{\color{blue} We should point out that, in general, an appropriately defined SEBEU in open-loop strategies, where each DM is restricted to selecting its control input only as a function of the time variable $t$, would not be an SEBEU in closed-loop strategies (in the sense of Definition~\ref{def:sebeu}). The reason for this is that an individual DM would face a stochastic optimal control problem given the strategies of the other DMs; therefore, an individual DM can typically reduce its cost by unilaterally deviating to a closed-loop strategy. However, in a deterministic setup, an SEBEU in open-loop strategies would also be an SEBEU in closed-loop strategies for the same reason. This is in contrast to the class of dynamic games where an open-loop Nash equilibrium is also a degenerate closed-loop Nash equilibrium; see for example \cite{fershtman1987identification}.}

\subsection{Examples and Intuition}\label{ss:ex}
\subsubsection{Example 1}
We present an example motivated by an autonomous demand response application, that is a single-stage finite-space game with no state variables. In such a game, a joint pure strategy $u=(u^1,\dots,u^N)$ determines the environment variable as $y=g(u,\xi)$, where $\xi$ is the random disturbance. As a result, each DM$^i$ incurs the cost $c^i(u^i,y)$.
A mixed strategy $\gamma^i$ for DM$^i$ is a probability mass function over DM$^i$'s set of pure strategies, whereas a subjective belief $\zeta^i$ for DM$^i$ is a probability mass function over the set of environment variables, i.e., $\zeta^i(y)$ is DM$^i$'s subjective probability for $y$.

In this setup, an SEBEU can be viewed as a profile of mixed strategies $\gamma=(\gamma^1,\dots,\gamma^N)$ and subjective beliefs $\zeta=(\zeta^1,\dots,\zeta^N)$ such that two conditions are met: (i) each DM$^i$'s mixed strategy $\gamma^i$ is optimal with respect to its subjective belief $\zeta^i$, i.e.,
$$\sum_{u^i\in\mathbb{U}^i} \gamma^i(u^i) \sum_{y\in\mathbb{Y}}  \zeta^i(y) c^i(u^i,y) \leq \sum_{u^i\in\mathbb{U}^i} \tilde{\gamma}^i(u^i) \sum_{y\in\mathbb{Y}} \zeta^i(y) c^i(u^i,y), \quad \forall  i\in[1,N], \tilde{\gamma}^i\in\Gamma^i$$
(ii) each DM$^i$'s subjective belief $\zeta^i$ is consistent with the objective probability $\zeta_{\gamma}$, i.e.,
$$\zeta^i(y)=\zeta_{\gamma}(y), \quad \forall  i\in[1,N], y\in\mathbb{Y}$$
where $\zeta_{\gamma}(y):=\sum_{(\bar{u},\bar{\xi}):g(\bar{u},\bar{\xi})=y} P[\bar{\xi}] \prod_{j\in[1,N]} \gamma^j(\bar{u}^j)$.

Let us now consider an autonomous demand response program where each DM$^i$ is an energy consumer who selects a consumption level $u^i\in\{0,1,2\}$. The unit price of energy (henceforth, simply ``the price''), which is the environment variable, equals the average consumption, i.e., $y=\frac{1}{N}\sum_{i\in[1,N]} u^i$. The net cost to each consumer equals the total price of the energy used  minus the utility derived by the consumer, which is given by $c^i(u^i,y)=u^i y - u^i$.

In this example, the consumers are at a pure-strategy SEBEU if and only if the average consumption equals one, i.e., $\frac{1}{N}\sum_{j\in[1,N]} u^j=1$, supported by the subjective beliefs satisfying $\zeta^1(1)=\cdots=\zeta^N(1)=1$.
At such an SEBEU, each consumer believes that the price equals to one with certainty (regardless of its consumption), and therefore the utility derived from any consumption would equal the total price of energy used by the consumer, i.e., $c^i(u^i,1)=0$, for all $u^i\in\{0,1,2\}$. Furthermore, because the average consumption equals one,  the true price is consistent with each consumer's subjective belief.
Note that a consumer at an SEBEU does not consider the possibility of its consumption changing the price, perhaps due to the lack of information on how the price is determined or the belief that the effect of its consumption on the price would be negligible (which would be true in a large scale program).

In contrast, a subjective equilibrium in the sense of \cite{Kalai:GEB:1995} would be a profile of mixed strategies $\gamma=(\gamma^1,\dots,\gamma^N)$ and subjective conditional beliefs $\zeta=\big((\zeta_{u^1}^1)_{u^1\in\mathbb{U}^1},\dots,$ $(\zeta_{u^N}^N)_{u^N\in\mathbb{U}^N}\big)$, where $\zeta_{u^i}^i$ is DM$^i$'s subjective probability mass function over the set of environment variables given $u^i$, such that two conditions are met: (i) each DM$^i$'s mixed strategy $\gamma^i$ is optimal with respect to its subjective belief $(\zeta_{u^i}^i)_{u^i\in\mathbb{U}^i}$, i.e.,
$$\sum_{u^i\in\mathbb{U}^i} \gamma^i(u^i) \sum_{y\in\mathbb{Y}}  \zeta_{u^i}^i(y) c^i(u^i,y) \leq \sum_{u^i\in\mathbb{U}^i} \tilde{\gamma}^i(u^i) \sum_{y\in\mathbb{Y}} \zeta_{u^i}^i(y) c^i(u^i,y), \quad \forall  i\in[1,N], \tilde{\gamma}^i\in\Gamma^i$$
(ii) each DM$^i$'s subjective conditional beliefs $(\zeta_{u^i}^i)_{u^i\in\mathbb{U}^i}$ are consistent with the objective conditional probabilities $(\zeta_{\gamma|u^i})_{u^i\in\mathbb{U}^i}$ as
\begin{equation}
\label{eq:SEconsistency}
\gamma^i(u^i) \zeta_{u^i}^i(y) = \gamma^i(u^i) \zeta_{\gamma|u^i}(y), \quad \forall  i\in[1,N], u^i\in\mathbb{U}^i, y\in\mathbb{Y}
\end{equation}
where $\zeta_{\gamma|u^i}(y):=\sum_{(\bar{u},\bar{\xi}):\bar{u}^i=u^i, g(\bar{u},\bar{\xi})=y} P[\bar{\xi}] \prod_{j\in[1,N],j\not=i} \gamma^j(\bar{u}^j)$.

In our example, the consumers would be at a pure-strategy subjective equilibrium in the sense of \cite{Kalai:GEB:1995} if and only if the average consumption is not higher than one, i.e., $\frac{1}{N}\sum_{j\in[1,N]} u^j\in\{0,1/N,\dots,1\}$. Such a subjective equilibrium can be supported, in particular, by any subjective conditional beliefs satisfying
$$\zeta_{u^i}^i\bigg(\frac{1}{N}\sum_{j\in[1,N]} u^j\bigg)=1 \ \ \textrm{and} \ \ E^{\zeta_{\tilde{u}^i}^i}[y]=\sum_{y\in\{0,1/N,\dots,2\}} y\zeta_{\tilde{u}^i}^i(y)\geq 1, \ \ \forall i\in[1,N], \tilde{u}^i \not= u^i.$$ Each consumer at such a subjective equilibrium holds the subjective conditional beliefs that (i) assigns probability one to the true average consumption conditioned on its true consumption, and (ii) leads to an expected price that is, in particular, not lower than one when conditioned on the alternative consumption levels. This ensures that its consumption is optimal with respect to each consumer's subjective conditional beliefs, which is consistent with the true price when conditioned on the consumer's true consumption. When conditioned on an alternative consumption level, a consumer's subjective belief need not match the objective reality and can result in higher expected price than the true price; in particular higher than one, in which case the consumer would have no incentive to switch to an alternative consumption level even when the price equals zero (i.e., all consumers are turned off).

When the consistency condition (\ref{eq:SEconsistency}) is replaced with $\zeta_{u^i}^i(y) = \zeta_{\gamma|u^i}(y)$, for all  $i\in[1,N]$, $u^i\in\mathbb{U}^i$, $y\in\mathbb{Y}$, then a Nash equilibrium is obtained. At a Nash equilibrium, the objective reality and the subjective beliefs of each DM are consistent when conditioned on every possible decision of the DM.
In our example, the consumers are at a Nash equilibrium if and only if every consumption level equals one, i.e., $u^1=\cdots=u^N=1$ or the average consumption equals $1-1/N$, i.e., $\frac{1}{N}\sum_{j\in[1,N]} u^j=1-1/N$. Note that a price $y \leq 1-2/N$ cannot arise at a NE because a consumer with $u^i=0$ (which necessarily exists to support such a price) can correctly predict that a unilateral switch to $\tilde{u}^i=1$ would increase the price to $\tilde{y}:=y+1/N$ and would decrease its net cost from zero to $\tilde{y}-1 <0$.
Similarly, a price $y \geq 1+1/N$ cannot arise at a NE because a consumer with $u^i>0$ (which necessarily exists to support such a price) can decrease its net cost from $u^i(y-1)>0$ to zero by unilaterally switching to $\tilde{u}^i=0$. The price $y=1$ arises at a Nash equilibrium if  $u^1=\cdots=u^N=1$; however, the price $y=1$ is not supported at a Nash equilibrium if there is a consumer with $u^i=2$ because such a consumer can unilaterally switch to $\tilde{u}^i=1$ and decrease its net cost from zero to $-1/N$. On the other hand, a Nash equilibrium arises whenever $y=1-1/N$ because (i) a consumer with $u^i=0$ cannot decrease its net cost by switching to $\tilde{u}^i\in\{1,2\}$ because the price would increase to $\tilde{y}\in\{1,1+1/N\}$, (ii) a consumer with $u^i=1$ who achieves the net cost $-1/N$ would increase its net cost to zero if it switches to $\tilde{u}^i\in\{0,2\}$, (iii) a consumer with $u^i=2$ who achieves the net cost $-2/N$ would still achieve the net cost $-2/N$ if it switches to $\tilde{u}^i=1$ or would achieve zero net cost if it switches to $\tilde{u}^i=0$.

This example illustrates that SEBEU is not only distinct from subjective equilibrium introduced in \cite{Kalai:GEB:1995} and Nash equilibrium but it can also be a more appropriate equilibrium concept in certain (large scale) applications where DMs cannot correctly discern the effects of their decisions on the environment variables due to a variety of reasons, in particular, due to lack of information. It should finally be pointed out that, at an SEBEU, a consumer who decides to model the price objectively by taking the effect of its own consumption on the price, can reduce its cost by at most $1/N$ by unilaterally switching to an alternative consumption level, i.e., an SEBEU is an approximate Nash equilibrium in a system with large number of consumers.

\subsubsection{Example 2}
We present another example where all state, control, disturbance, and environment variables take scalar real values. Each DM$^i$ has its own linear state dynamics
$$x_{t+1}^i  = a x_t^i +b u_t^i  + w_t^i, \qquad t\in[0,T)$$
and its own quadratic cost
$$ \sum_{t\in[0,T)} E\big[q(x_t^i)^2 + r(u_t^i)^2 + 2y_t  u_t^i\big]  + qE\big[(x_T^i)^2\big]$$
where $a$, $b$, $q>0$, $r>0$ are scalar constants and $T\in\mathbb{N}$.
The environment variables are generated by the linear dynamics
\begin{align*}
y_t = &  \frac{1}{N} \sum_{i\in[1,N]} u_t^i  +\xi_t
\end{align*}
We assume $x_0^1,\dots,x_0^N,w_0^1,\dots,w_0^N,\dots,w_{T-1}^1,\dots,w_{T-1}^N,\xi_0,\dots,\xi_{T-1}\in\mathbb{R}$ are mutually independent and zero-mean Gaussian where
$\textrm{cov}[\xi_0]>0$.

If DM$^i$ views the environment variables as an independent exogenous random sequence $(z_t)_{t\in[0,T)}$ (the random variables $z_0,\dots,z_{T-1}$ need not be independent), then its unique optimal strategy can be derived by dynamic programming as, for $k\in[0,T)$,
\begin{align*}
s_k^i(I_k^i,Z_{k-1}) = & f_k x_k^i +  \sum_{t\in[k,T)} g_{k,t} E[z_t|Z_{k-1}]
\end{align*}
where $Z_{-1}=\emptyset$, $Z_{k-1}=(z_0,\dots,z_{k-1})$, $k\in[1,T)$, and
\begin{align*}
f_k = & -\frac{b m_{k+1} a}{r+ m_{k+1} b^2}, \quad
g_{k,k} =  - \frac{1}{r+ m_{k+1} b^2}, \quad
g_{k,t} =  -\frac{b \phi_{k,t-1} f_t}{r+ m_{k+1} b^2} , \ \ t \in [k+1,T)     \\
m_k = & q +  \frac{r m_{k+1} a^2}{r+ m_{k+1} b^2}   , \ \ m_T=q \\
\phi_{k,k} =& 1, \quad \phi_{k,t} =  \prod_{n\in[k+1,t]} (a + b f_{n}), \ \ t\in[k+1,T).
\end{align*}

In view of the above, it follows that a joint strategy $s=(s^1,\dots,s^N)\in\mathbb{S}$ is an SEBEU if and only if, for $k\in[0,T)$,
\begin{align*}
s_k^i(I_k^i,Y_{k-1}) = & f_k x_k^i +  \sum_{t\in[k,T)} g_{k,t} E[y_t|Y_{k-1}]
\end{align*}
where $(y_0,\dots,y_{T-1})$ satisfies
\begin{align}
y_k = &  f_k  \bar{x}_k +  \sum_{t\in[k,T)} g_{k,t} E[y_t|Y_{k-1}]   +\xi_k \label{eq:yk} \\
\bar{x}_{k+1}  = & (a + bf_k) \bar{x}_k + b \sum_{t\in[k,T)} g_{k,t} E[y_t|Y_{k-1}]   + \bar{w}_k \label{eq:barxkp1}
\end{align}
and $\bar{x}_0=\frac{1}{N} \sum_{i\in[1,N]} x_0^i$, $\bar{w}_k=\frac{1}{N} \sum_{i\in[1,N]} w_k^i$.

Now, let us consider the special case of $T=2$. From (\ref{eq:yk})-(\ref{eq:barxkp1}), we obtain
\begin{align*}
\left[\begin{array}{cc}
1 - g_{0,0} & -g_{0,1}  \\
-f_1 b g_{0,0} &  1 - f_1 b g_{0,1} -  g_{1,1}
\end{array} \right] \left[\begin{array}{c}
E[y_0] \\ E[y_1]
\end{array} \right]=
0.
\end{align*}
Since $(1 - g_{0,0})(1 - f_1 b g_{0,1} -  g_{1,1})-f_1 b g_{0,0}g_{0,1}>0$, we must have $E[y_0]=E[y_1]=0$.
From (\ref{eq:yk})-(\ref{eq:barxkp1}), we also obtain
\begin{align*}
E[y_1|y_0] = &  \frac{f_1  E[\bar{x}_1|y_0]}{1-g_{1,1}}
\end{align*}
where $1-g_{1,1}>0$ and
\begin{align*}
E[\bar{x}_1|y_0]  = & \frac{(a + bf_0) \textrm{cov}[\bar{x}_0] f_0}{\textrm{cov}[\bar{x}_0]f_0^2+\textrm{cov}[\xi_0]} y_0.
\end{align*}
Therefore, the unique SEBEU is given by, $i\in[1,N]$,
\begin{align}
s_0^i(I_0^i) = & f_0 x_0^i \label{eq:ex2se0} \\
s_1^i(I_1^i,Y_0) = & f_1 x_1^i +    k_1^N y_0 \label{eq:ex2se1}
\end{align}
where
\begin{align*}
k_1^N := &    \frac{  - 1 }{1+ r+ q b^2} \frac{ f_1 (a + bf_0) \textrm{cov}[\bar{x}_0] f_0 }{\textrm{cov}[\bar{x}_0]f_0^2+\textrm{cov}[\xi_0]}.
\end{align*}
The state feedback gains $f_0,f_1$ are independent of the number of DMs; furthermore, $\lim_{N\rightarrow\infty} k_1^N=0$ due to $\textrm{cov}[\bar{x}_0]=\frac{1}{N}\textrm{cov}[x_0^1]$. In the case of $T=2$, obtaining the set of Nash equilibria is not as straightforward; however, it is possible to show that, for any $\epsilon>0$, the SEBEU (\ref{eq:ex2se0})-(\ref{eq:ex2se1}) constitutes an $\epsilon-$Nash equilibrium if the number of DMs is sufficiently large.
To show this, we obtain DM$^i$'s optimal response $\tilde{s}^{i,N}\in\mathbb{S}^i$ to the strategies (\ref{eq:ex2se0})-(\ref{eq:ex2se1}) as follows.
\begin{align*}
\tilde{V}_1^i(I_1^i,y_0) := & \min_{u_1^i} E[q(x_1^i)^2+r(u_1^i)^2 + 2u_1^iy_1 + q(x_2^i)^2 | I_1^i,y_0] \\
                  = & \min_{u_1^i} q(x_1^i)^2+r(u_1^i)^2   +2(u_1^i)^2/N +2u_1^i f_1 (a+bf_0)   E[\bar{x}_0^{-i}|I_1^i,y_0] \\
                    &  \qquad + 2u_1^ik_1 y_0(N-1)/N  +  qE \big[ (a x_1^i + b u_1^i + w_1^i)^2 | I_1^i,y_0 \big]
\end{align*}
where $\bar{x}_0^{-i}:=\frac{1}{N} \sum_{j\not=i}  x_0^j$. DM$^i$'s optimal response at stage $t=1$ is
$$\tilde{s}_1^{i,N}(I_1^i,y_0)=-\frac{qbax_1^i+ f_1 (a+bf_0)   E[\bar{x}_0^{-i}|I_1^i,y_0] +  k_1 y_0(N-1)/N}{r+qb^2+2/N}. $$
From $y_0= \frac{u_0^i}{N} + f_0 \bar{x}_0^{-i} + \xi_0$, we have
$$E[\bar{x}_0^{-i}|I_1^i,y_0]=\frac{f_0\textrm{cov}[\bar{x}_0^{-i}]}{f_0^2\textrm{cov}[\bar{x}_0^{-i}]+\textrm{cov}[\xi_0]}(y_0-u_0^i/N).$$
Therefore,
$$\tilde{s}_1^{i,N}(I_1^i,y_0)=\tilde{f}_1^N x_1^i +\tilde{k}_1^N y_0 +\tilde{n}_1^N u_0^i $$
where
\begin{align*}
\tilde{f}_1^N  := & -\frac{qba}{r+qb^2+2/N} \\
\tilde{k}_1^N     := & -\frac{f_1 (a+bf_0)    }{r+qb^2+2/N} \frac{f_0\textrm{cov}[\bar{x}_0^{-i}]}{f_0^2\textrm{cov}[\bar{x}_0^{-i}]+\textrm{cov}[\xi_0]}
-\frac{ k_1(N-1)/N }{r+qb^2+2/N}
\\
\tilde{n}_1^N     := & \frac{1}{N} \frac{f_1 (a+bf_0)    }{r+qb^2+2/N} \frac{f_0\textrm{cov}[\bar{x}_0^{-i}]}{f_0^2\textrm{cov}[\bar{x}_0^{-i}]+\textrm{cov}[\xi_0]}.
\end{align*}
%This results in
%\begin{align*}
%\tilde{V}_1^i(I_1^i,y_0) = &  q(x_1^i)^2+r(\tilde{f}_1^N x_1^i +\tilde{k}_1^N y_0 +\tilde{n}_1^N u_0^i)^2   \\
%                    & +2\frac{1}{N} (\tilde{f}_1^N x_1^i +\tilde{k}_1^N y_0 +\tilde{n}_1^N u_0^i)^2 \\
%                    & +2 (\tilde{f}_1^N x_1^i +\tilde{k}_1^N y_0 +\tilde{n}_1^N u_0^i) f_1 (a+bf_0)   \frac{f_0\textrm{cov}(\bar{x}_0^{-i})}{f_0^2\textrm{cov}(\bar{x}_0^{-i})+\textrm{cov}(\xi_0)}(y_0-u_0^i/N) \\
%                    & +  2 (\tilde{f}_1^N x_1^i +\tilde{k}_1^N y_0 +\tilde{n}_1^N u_0^i) k_1 y_0 (N-1) /N \\
%                    & +  qE ((a x_1^i + b (\tilde{f}_1^N x_1^i +\tilde{k}_1^N y_0 +\tilde{n}_1^N u_0^i) + w_1^i)^2 |I_1^i,y_0)
%\end{align*}
Let
\begin{align*}
\tilde{V}_0^i(x_0^i) := & \min_{u_0^i} E\big[q(x_0^i)^2+r(u_0^i)^2   +2u_0^i y_0+  \tilde{V}_1^i (I_1^i,y_0) | x_0^i)\big].
\end{align*}
%DM$^i$'s optimal response at stage $t=0$ satisfies
%\begin{align*}
%0 = & 2ru_0^i   +4 \frac{u_0^i}{N} +  2qb(ax_0^i+bu_0^i) \\
%    & + (2r+4/N)(\tilde{f}_1^N (ax_0^i + bu_0^i) +\tilde{k}_1^N u_0^i / N +\tilde{n}_1^N u_0^i)(\tilde{f}_1^N b +\tilde{k}_1^N/N  +\tilde{n}_1^N) \\
%    & +  2 (\tilde{f}_1^N  b +\tilde{k}_1^N /N +\tilde{n}_1^N ) ((N-1)/N) k_1 u_0^i/N \\
%    & +  2 (\tilde{f}_1^N (a x_0^i + bu_0^i ) +\tilde{k}_1^N u_0^i/N +\tilde{n}_1^N u_0^i) ((N-1)/N)k_1 /N \\
%                    & +  2q((a + b\tilde{f}_1^N) (ax_0^i + bu_0^i) + b (\tilde{k}_1^N u_0^i/N +\tilde{n}_1^N u_0^i))
%                           ((a + b\tilde{f}_1^N) b + b (\tilde{k}_1^N /N +\tilde{n}_1^N ))
%\end{align*}
DM$^i$'s optimal response at stage $t=0$ is obtained as
$$\tilde{s}_0^{i,N}(x_0^i)=\tilde{f}_0^N x_0^i$$
where
\begin{align*}
\tilde{f}_0^N := & -n_{\tilde{f}_0^N} / d_{\tilde{f}_0^N} \\
n_{\tilde{f}_0^N} := &   qba + (r+2/N)\tilde{f}_1^N a (\tilde{f}_1^N b +\tilde{k}_1^N/N  +\tilde{n}_1^N) +  2 \tilde{f}_1^N a    ((N-1)/N) k_1 /N \\
                    & +  q(a + b\tilde{f}_1^N) ab   (a + b\tilde{f}_1^N  +  \tilde{k}_1^N /N +\tilde{n}_1^N )   \\
d_{\tilde{f}_0^N} := & r   + 2 /{N} +  qb^2 + (r+2/N)(\tilde{f}_1^N b +\tilde{k}_1^N / N +\tilde{n}_1^N )^2 \\
    & +  2 (\tilde{f}_1^N  b +\tilde{k}_1^N /N +\tilde{n}_1^N ) ((N-1)/N) k_1 /N   +  qb^2(a + b\tilde{f}_1^N + \tilde{k}_1^N /N +\tilde{n}_1^N )^2.
\end{align*}
It is straightforward to see that $\lim_{N\rightarrow\infty}(\tilde{f}_0^N,\tilde{f}_1^N,\tilde{k}_1^N,\tilde{n}_1^N)= (f_0,f_1,0,0)$.
This implies that, for any $\epsilon>0$, the SEBEU (\ref{eq:ex2se0})-(\ref{eq:ex2se1}) is an $\epsilon-$Nash equilibrium if the number of DMs is sufficiently large.

We should also remark that the SEBEU strategies (\ref{eq:ex2se0})-(\ref{eq:ex2se1}) are {\it always} linear mappings of each DM$^i$'s information variables whereas it is not clear whether the same holds for Nash equilibria for such linear quadratic Gaussian models \cite[Chapter 7]{basols99}.

\section{Existence of SEBEU in mixed strategies}
In this section, we consider the general setup introduced in Section~\ref{se:model}. Our goal is to establish the existence of SEBEU in mixed strategies when all state, control, disturbance, and environment variables belong to compact metric spaces.

\begin{assumption}\label{as:cm}
$\textrm{}$

\begin{enumerate}
\item Each of the sets $\mathbb{X}^0,\dots,\mathbb{X}^N,\mathbb{U}^1,\dots,\mathbb{U}^N,\mathbb{W}^0,\dots,\mathbb{W}^N,\mathbb{Y},\Xi$
 is a compact metric space
\item For all $(i,t)\in[0,N]\times[0,T]$, $c_t^i$, $f_t^i$, $g_t$ are continuous (if defined)
\item For all $i\in[1,N]$, $\mathbb{S}^i=\times_{t\in[0,T)} \mathbb{S}_t^i$ where, for each $t\in[0,T)$, $\mathbb{S}_t^i$ is a compact subset of the space of continuous functions from $\mathbb{I}_t^i \times \mathbb{Y}^{[0,t)}$ to $\mathbb{U}^i$ equipped with the supremum norm
%    \footnote{The uniform metric for the space $\mathcal{C}$ of continuous functions from $\mathbb{A}$ to $\mathbb{B}$, where $\mathbb{A}$ and $\mathbb{B}$
%    are metric spaces, is defined by  $m_{\mathcal{C}}(f,g) = \sup_{a\in\mathbb{A}} m_{\mathbb{B}} (f(a),g(a))$ where $m_{\mathbb{B}}$ is the metric on
%    $\mathbb{B}$.}
\item For some $\bar{c}<\infty$ and $\bar{\beta}<1$, $\sup_{(i,t,x^i,u^i,y)\in[1,N]\times[0,T]\times\mathbb{X}^i\times\mathbb{U}^i\times\mathbb{Y}} |c_t^i(x^i,u^i,y)|<\bar{c}$ and
$\sup_{i\in[1,N]} \beta^i < \bar{\beta}$ .
\end{enumerate}
\end{assumption}

Condition (1) of Assumption~\ref{as:cm} is satisfied by finite sets equipped with the discrete metric, in which case Condition (2) is satisfied by all functions $c_t^i$, $f_t^i$, $g_t$ and Condition (3) is satisfied if $\mathbb{S}_t^i$ is the set of all mappings from $\mathbb{I}_t^i \times \mathbb{Y}^{[0,t)}$ to $\mathbb{U}^i$.
%Note that the functions $c_t^i$, $f_t^i$,  $g_t$ satisfying Condition (2) are uniformly continuous  (in their compact domains of definition).
For more general spaces, Condition (3) is satisfied, for example, by the space of Lipschitz continuous functions from $\mathbb{I}_t^i \times \mathbb{Y}^{[0,t)}$ to $\mathbb{U}^i$, with a uniform Lipschitz constant, equipped with the supremum norm.
Note that, for all $t\in[0,T)$, the compactness of
$\mathbb{I}_t^i \times \mathbb{Y}^{[0,t)}$  and
$\mathbb{S}_t^i$ imply that the mapping in $\mathbb{S}_t^i$ are uniformly equi-continuous
    % \footnote{\label{fn:eqcondef}
    % A collection $\mathcal{F}$ of continuous functions from $\mathbb{A}$ to $\mathbb{B}$, where $\mathbb{A}$ and $\mathbb{B}$ are metric spaces, is called
    % uniformly equi-continuous if, for all $\epsilon>0$, there exists $\delta>0$ such that, for all $a\in\mathbb{A}$ and $f\in\mathcal{F}$,
    % $$m_{\mathbb{A}}(a  ,  \bar{a}) < \delta \qquad\Rightarrow\quad m_{\mathbb{B}}(f(a)  ,  f(\bar{a}) < \epsilon.$$}
by the Arzela-Ascoli Theorem \cite{Dud02}. Condition (4) implies that any DM$^i$'s long-term cost is uniformly bounded, i.e., $\sup_{(i,s^i,z)\in[1,N]\times\mathbb{S}^i\times\mathbb{Y}^{[0,T)}} |J^i(s^i,z)| <\infty$, which will be needed in the case $T=\infty$.

We now turn to the question of existence of SEBEU in mixed strategies.
If DM$^i$ faces an independent exogenous distribution $\zeta\in\mathcal{P}(\mathbb{Y}^{[0,T)})$ for the environment variables, DM$^i$'s set of best replies to $\zeta$ in mixed strategies are defined as
$$BR^i(\zeta):=\{\gamma^i\in\Gamma^i: E^{(\gamma^i,\zeta)} [J^i] \leq E^{(\tilde{\gamma}^i,\zeta)} [J^i], \ \forall \tilde{\gamma}^i\in\Gamma^i \}.$$
We now introduce the correspondence $\overline{BR}=(\overline{BR}^1,\dots,\overline{BR}^N):\Gamma\rightrightarrows\Gamma$ defined by
$$\overline{BR}^i(\gamma)=BR^i(\zeta_{\gamma}), \quad \forall i\in[1,N].$$ Clearly, a mixed strategy profile $\gamma\in\Gamma$ is an SEBEU if and only if $\gamma$ is a fixed point of $\overline{BR}$.\\

\begin{theorem}
\label{thm:cm}
Consider the compact metric space model of this section, under Assumption~\ref{as:cm}.
(i) There exists an SEBEU in mixed strategies.
%If DMs are identical, there exists a \textrm{symmetric} SEBEU in mixed strategies $\gamma\in\Gamma$ such that $\gamma^1=\cdots=\gamma^N$.
(ii) If there are sets of identical\footnote{DM$^i$ and DM$^j$ are considered identical if $BR^i(\zeta)=BR^j(\zeta)$, for all $\zeta\in\mathcal{P}(\mathbb{Y}^{[0,T)})$.} DMs, there exists a symmetric SEBEU in mixed strategies $\gamma\in\Gamma$ such that identical DMs have identical strategies at  $\gamma$, i.e., $\gamma^i=\gamma^j$, if DM$^i$ and DM$^j$ are identical.\\
\end{theorem}

A proof for Theorem~\ref{thm:cm} based on Fan-Glicksberg's fixed point theorem \cite{Fan1952,Glicksberg1952} is provided in Appendix~\ref{appx:ex}.

Due to the particular product-space (across time-stages) construction of the individual set of pure strategies and the mixed strategies being probability distributions over such sets, it may be useful to motivate the practicality of the existence result above. Toward this end, we recall that
a \textit{behavioral strategy} for DM$^i$ is a collection $b^i=(b_t^i)_{t\in[0,T)}$ of functions such that  $b_t^i:\mathbb{I}_t^i \times\mathbb{Y}^{[0,T)} \to \mathcal{P}(\mathbb{U}^i)$, for all $t\in[0,T)$. If DM$^i$ employs a behavior strategy $b^i$, then DM$^i$ makes an independent draw $u^i\in\mathbb{U}^i$ at each time $t\in[0,T)$ according to $b_t^i(I_t^i,Y_{t-1})$ and implements $u_t^i=u^i$.   It is well known that, in extensive games with imperfect information such as the games considered in this paper, if the perfect recall condition holds, i.e., each DM$^i$ never forgets any information, then every mixed strategy has an equivalent behavior strategy; see \cite{kuhn1953extensive,aumann1961mixed}.

More precisely, given a mixed strategy $\gamma^i\in\Gamma^i$ for DM$^i$, an (essentially unique) equivalent behavior strategy can be defined as: for all $t\in[0,T)$, $(I_t^i,Y_{t-1})\in\mathbb{I}_t^i\times\mathbb{Y}^{[0,T)}$, $B^i\in\mathcal{B}(\mathbb{U}^i)$,
\begin{align}
\label{eq:bs1}
b_t^i(I_t^i,Y_{t-1})(B^i) & := P^{\gamma^i} \big[s_t^i(I_t^i,Y_{t-1})\in B^i \ | \ s^i\in R_{t}^i(I_{t}^i,Y_{t-1})\big]
\end{align}
where $R_t^i$ denotes the set of DM$^i$'s pure strategies consistent with $(I_{t}^i,Y_{t-1})$, i.e.,
$$R_t^i(I_t^i,Y_{t-1}):=\{s^i\in\mathbb{S}^i:s_k^i(I_k^i,Y_{k-1})=u_k^i, \ k\in[0,t)\}$$
provided that $P^{\gamma^i}[R_{t}^i(I_{t}^i,Y_{t-1})]>0$.  If $P^{\gamma^i}[R_{t}^i(I_{t}^i,Y_{t-1})]=0$, the history $(I_t^i,Y_{t-1})$ arises with zero probability under $\gamma^i$ (regardless of the strategies employed by the other DMs) in which case the definition of $b_t^i(I_t^i,Y_{t-1})$ is arbitrary.

It is straightforward to verify that a mixed strategy profile $\gamma=(\gamma^1,\dots,\gamma^N)\in\Gamma$ and its equivalent behavior strategy $b=(b^1,\dots,b^N)$ defined by (\ref{eq:bs1}) induce the same probability distribution over the set $\mathbb{I}_T^1\times\cdots\times\mathbb{I}_T^N\times\mathbb{Y}^{[0,T)}$ of histories.
Similarly, $(s^i,Z_{T-1})\in\mathbb{S}^i\times\mathbb{Y}^{[0,T)}$ and its equivalent $(b^i,Z_{T-1})$ induce the same probability distribution over the set $\mathbb{I}_T^i$ of histories. This leads to the following corollary to Theorem~\ref{thm:cm}.\\

\begin{corollary}
\label{cor:cm}
Consider the compact metric space model of this section under Assumption~\ref{as:cm}.
(i) There exists an SEBEU in behavior strategies.
(ii) If there are sets of identical DMs, there exists a symmetric SEBEU in behavior strategies $b$ such that identical DMs have identical strategies at  $b$.\\
\end{corollary}

%%%%%%%%%%%%%%%%%%%%%%%%%%%%%%%%%%%%%%%%%%%%%%%%%%%%%%%%%%%%%%%%%%%%%%%%%%%%%%%%
\section{Large games with small players}
\label{se:anongame}
In this section, we consider a class of games in which a large number of ``small'' players interact within the compact metric space model of the previous section; see \cite{fudenberg1998nonanonymous,BLONSKI2000225,adlakha2015equilibria} for the case of finite set of players and \cite{mas1984theorem, schmeidler1973equilibrium, jovanovic1988anonymous} for the case of non-atomic set of players.
%\ga{We do not require DMs to be identical or anonymous; however, for technical reasons, we assume that they share a common set of individual states and a common set of individual controls, i.e.,  $\mathbb{X}^1=\cdots=\mathbb{X}^N$ and $\mathbb{U}^1=\cdots=\mathbb{U}^N$.}
%
%In this section, we consider a special class of games, called ``anonymous games'' (see \cite{BLONSKI2000225,adlakha2015equilibria} for the case of finite set of players and \cite{mas1984theorem, schmeidler1973equilibrium, jovanovic1988anonymous} for the case of non-atomic set of players), within the compact metric space model of the previous section.
%
As the notion of an SEBEU is different from that of a Nash equilibrium, the long-term cost of each DM at an SEBEU need not be optimal given the strategies of the other DMs. However, in large games where each DM's impact on the environment variables is small, the optimality gap arising at an SEBEU vanishes asymptotically  as the number of DMs grows unboundedly.
%To introduce the notion of an anonymous game, we consider the case where the set of states and control inputs available to each DM is identical, i.e., %$(\mathbb{X}^i\times\mathbb{U}^i)=(\mathbb{X}^1\times\mathbb{U}^1)$, for all $i\in[1,N]$.
In this section, we consider games in which the environment variables $(y_t)_{t\in[0,T)}$ depend on the state-control variables $((x_t,u_t))_{t\in[0,T)}$ through some aggregate statistics $(d_t)_{t\in[0,T)}$ which take values in some compact metric space $\mathbb{D}$. For example, when $\mathbb{X}^1=\cdots=\mathbb{X}^N \subset \mathbb{R}^n$ and $\mathbb{U}^1=\cdots=\mathbb{U}^N \subset \mathbb{R}^m$ for some $n,m\in\mathbb{N}$, the aggregate statistics can be
\begin{equation}
\label{eq:avgfun}
d_t= \frac{1}{N}\sum_{i\in[1,N]} \begin{bmatrix} x_t^i \\ u_t^i \end{bmatrix}, \quad t\in[0,T).
\end{equation}

The aggregate statistics $(d_t)_{t\in[0,T)}$ can be different from the simple averages in (\ref{eq:avgfun}) as long as the following hold: (i) regardless of the number of DMs in the system, DMs can cause only small changes in the aggregate statistics if the maximum change across all DMs is small (\ref{eq:negdm1}), and (ii)  the individual impact of any DM$^i$ on the aggregate statistics vanishes as the number of DMs grows (\ref{eq:negdm2}).
In this section, the equations (\ref{eq:p})-(\ref{eq:y}) generating the environment variables take the form
\begin{align}
\label{eq:p2}
y_t =  & g_t(x_t^0, d_t,\xi_t), \qquad t\in[0,T) \\
\label{eq:y2}
x_{t+1}^0 = & f_t^0(x_t^0, d_t, w_t^0), \qquad t\in[0,T)
\end{align}
where $g_t:\mathbb{X}^0\times\mathbb{D}\times\Xi\to\mathbb{Y}$, $f_t^0:\mathbb{X}^0\times\mathbb{D}\times\mathbb{W}^0\to\mathbb{X}^0$, $t\in[0,T)$.

To state our results precisely, we consider a sequence of $(\mathbb{G}^N)_{N\in\mathbb{N}}$ of games, where the number of DMs in $\mathbb{G}^N$ is $N\in\mathbb{N}$.
For ease of notation, we assume that all games in $(\mathbb{G}^N)_{N\in\mathbb{N}}$ share the common data $(T,\mathbb{Y},\mathbb{X}^0,\mathbb{D},\mathbb{W}^0,\Xi,(f_t^0)_{t\in[0,T)},(g_t)_{t\in[0,T)})$; however, the definition of the aggregate statistics $(d_t)_{t\in[0,T)}$ necessarily depends on the number $N$ of DMs as $d_t$ maps $\mathbb{X}^1\times\cdots\mathbb{X}^N\times\mathbb{U}^1\times\cdots\mathbb{U}^N$ to $\mathbb{D}$ in $\mathbb{G}^N$; hence we use the notation $d_t^N$ for the aggregate statistics at time $t\in[0,T)$ in the game $\mathbb{G}^N$, $N\in\mathbb{N}$,
as well as the notation $d_t^N(x_t,u_t)$ to make the dependence of $d_t^N$ on $(x_t,u_t)$ explicit.
Also, the data specific to any DM$^i$, $i\in\mathbb{N}$, that is, $(\mathbb{X}^i,\mathbb{U}^i,\mathbb{W}^i,(f_t^i)_{t\in[0,T)}$, $(c_t^i)_{t\in[0,T]},\beta^i)_{i\in[1,N]}$ as well as the distribution of DM$^i$'s local random parameters $(x_0^i,(w_t^i)_{t\in[0,T)})$ (which will be assumed independent across all DMs) remain the same in all games $\mathbb{G}^i,\mathbb{G}^{i+1},\dots$ containing DM$^i$. In other words, the game $\mathbb{G}^{N+1}$ is obtained by adding a new DM, that is DM$^{N+1}$, to the game $\mathbb{G}^N$ and modifying the definition of the aggregate statistics $(d_t^N)_{t\in[0,T)}$ so that each $d_t^{N+1}$, $t\in[0,T)$, maps $\mathbb{X}^1\times\cdots\mathbb{X}^{N+1}\times\mathbb{U}^1\times\cdots\mathbb{U}^{N+1}$ to $\mathbb{D}$.
We also use the notation $\bar{J}^{i,N}$ to denote $\bar{J}^i$ in $\mathbb{G}^N$ as defined in (\ref{eq:ltc1}) ($J^i$ defined in (\ref{eq:ltc2}) is independent of $N$, i.e., it is the same function  in all games $\mathbb{G}^i,\mathbb{G}^{i+1},\dots$); however we will often suppress the dependence of the various terms on $N$.    \\

\begin{assumption}
\label{as:anon1}
$\textrm{}$

\begin{enumerate}
\item Assumption~\ref{as:cm} holds  in all games $(\mathbb{G}^N)_{N\in\mathbb{N}}$ with the following strengthened conditions:
\begin{enumerate}
\item For all $i\in\mathbb{N}$, $\mathbb{X}^i=\mathbb{X}^1$,  $\mathbb{U}^i=\mathbb{U}^1$, $\mathbb{W}^i=\mathbb{W}^1$, and $\mathbb{S}^i=\mathbb{S}^1$
\item For all $t\in [0,T]$, $(c_t^i)_{i\in\mathbb{N}}$, $(f_t^i)_{i\in\mathbb{N}}$ are uniformly equi-continuous (if defined)
\item For some $\bar{c}<\infty$ and $\bar{\beta}<1$, $$\sup_{(i,t,x^i,u^i,y)\in\mathbb{N}\times[0,T]\times\mathbb{X}^i\times\mathbb{U}^i\times\mathbb{Y}} |c_t^i(x^i,u^i,y)|<\bar{c} \quad\textrm{and}\quad
\sup_{i\in\mathbb{N}} \beta^i < \bar{\beta}$$
\end{enumerate}
\item  For all $t\in[0,T)$ and $\epsilon>0$, there exists $\delta_{\epsilon}>0$ such that, in all games $(\mathbb{G}^N)_{N\in\mathbb{N}}$,
\begin{equation}
\label{eq:negdm1}
\sup_{i\in[1,N]} m_{\mathbb{X}^i\times\mathbb{U}^i}((x^i,u^i),(\tilde{x}^i,\tilde{u}^i)) < \delta_{\epsilon} \quad  \Rightarrow\quad m_{\mathbb{D}} (d_t^N(x,u),d_t^N(\tilde{x},\tilde{u}))<\epsilon
\end{equation}
and
\begin{equation}
\label{eq:negdm2}
\lim_{N\rightarrow\infty} \sup_{i\in[1,N]} \sup_{\scriptsize
\begin{array}{c}(x,u,\tilde{x},\tilde{u})\in\mathbb{X}\times\mathbb{U}\times\mathbb{X}\times\mathbb{U}\\:(x^{-i},u^{-i})=(\tilde{x}^{-i},\tilde{u}^{-i})\end{array}} m_{\mathbb{D}} (d_t^N(x,u),d_t^N(\tilde{x},\tilde{u}))=0
\end{equation}
\item The random variables $\omega^0,\omega^1,\omega^2,\dots$ defined as $$\omega^0:=(x_0^0,(w_t^0)_{t\in[0,T)},(\xi_t,)_{t\in[0,T)}), \ \omega^i:=(x_0^i,(w_t^i)_{t\in[0,T)}), \ i\in\mathbb{N}$$
    are mutually independent ($\Omega^0:=\mathbb{X}^0 \times (\mathbb{W}^0)^{[0,T)}\times(\Xi)^{[0,T)}$, $\Omega^i:=\mathbb{X}^i\times(\mathbb{W}^i)^{[0,T)}$, $i\in\mathbb{N}$).\\
\end{enumerate}
\end{assumption}

Condition~1 of Assumption~\ref{as:anon1} is to ensure the existence of an SEBEU in every game $\mathbb{G}^N$, $N\in\mathbb{N}$;
the strengthened conditions are needed to prove the results of this section where the number of DMs grows without bound, i.e., $N\rightarrow\infty$.
In particular, Condition~1(a) is to ensure that each DM$^i$'s set $\mathbb{I}_T^i \times \mathbb{Y}^{[0,T)}$ of observation histories is the same and so is each DM$^i$'s set $\mathbb{S}^i$ of pure strategies (the condition $\mathbb{W}^i=\mathbb{W}^1$, for all $i\in\mathbb{N}$, is for notational convenience and can be relaxed), whereas Condition~1(b)-(c) are to ensure that the long-term cost functions of DMs can be approximated by a finite set of functions with arbitrary precision.
Condition~2 is to ensure that (i) the change in the environment variables caused by a uniformly small change in all state-control variables is small regardless of the number of DMs in the system (\ref{eq:negdm1}),  and (ii) the change in the environment variables caused by a single DM's deviation to an alternative strategy is arbitrarily small in games with sufficiently large number of DMs (\ref{eq:negdm2}). For example, the simple averages in (\ref{eq:avgfun}) satisfy (\ref{eq:negdm1})-(\ref{eq:negdm2}). Intuitively, in a large game, the influence of an individual DM on the environment variables hence on the long-term cost of the other DMs would be negligible.  Note that we do not require DMs to be identical or anonymous.
Condition~3 is to simplify the proof of the next result, which states that an SEBEU well approximates a Nash equilibrium in large games, i.e., the optimality gap for each individual DM at an SEBEU vanishes as the number of DMs grow unboundedly.

\begin{theorem}
\label{th:exanon1}
Consider a sequence $(\mathbb{G}^N)_{N\in\mathbb{N}}$ of games under Assumption~\ref{as:anon1}.
Given $\epsilon>0$, there exist $\bar{N}_{\epsilon}\in\mathbb{N}$  such that, if $N\geq \bar{N}_{\epsilon}$,
then any SEBEU $\gamma^{*,N}=(\gamma^{1,N},\dots,\gamma^{N,N})$ in $\mathbb{G}^N$ is also an $\epsilon-$Nash equilibrium in $\mathbb{G}^N$, i.e.,
\begin{equation}
E^{\gamma^{*,N}} [\bar{J}^{i,N}] \leq E^{(\tilde{\gamma}^i,\gamma^{-i,N})}[\bar{J}^{i,N}] + \epsilon, \quad \forall \ i\in[1,N], \ \tilde{\gamma}^i\in\Gamma^i.
\label{eq:anonnash}
\end{equation}
\end{theorem}

\textbf{Proof of Theorem~\ref{th:exanon1}}
It suffices to show that
\begin{equation}
\label{eq:anonnashapprox}
\lim_{N\rightarrow\infty} \sup_{i\in[1,N],\tilde{\gamma}^i\in\Gamma^i,\gamma\in\Gamma} |E^{(\tilde{\gamma}^i,\zeta_{\gamma})} [ J^i ] - E^{(\tilde{\gamma}^i,\gamma^{-i})} [\bar{J}^{i,N}]| =0.
\end{equation}
When $T=\infty$, due to Condition~1(c) of Assumption~\ref{as:anon1}, we have
\begin{equation}
\label{eq:costtruncinf} \lim_{\bar{T}\rightarrow\infty} \sup_{i\in\mathbb{N},(I_T^i,Y_{T-1}) \in \mathbb{I}_T^i \times \mathbb{Y}^{[0,T)}} \Bigg|\sum_{t\in[\bar{T},\infty)} (\beta^i)^t c_t^i\big(x_t^i,u_t^i,z_t\big)  \Bigg| = 0.
\end{equation}
Therefore, it suffices to show (\ref{eq:anonnashapprox}) for the case of $T\in\mathbb{N}$.

For $i\in[1,N]$, $\tilde{\gamma}^i\in\Gamma^i$, $\gamma\in\Gamma$, let $\mu^{i,(\tilde{\gamma}^i,\zeta_{\gamma})}$ and $\mu^{i,(\tilde{\gamma}^i,\gamma^{-i})}$ denote the probability measures induced by $(\tilde{\gamma}^i,\zeta_{\gamma})$ and $(\tilde{\gamma}^i,\gamma^{-i})$, respectively, on  $\mathcal{B}(\mathbb{I}_T^i \times \mathbb{Y}^{[0,T)})$, where the environment variables generated by $\zeta_{\gamma}\in\mathcal{P}(\mathbb{Y}^{[0,T)})$ is assumed to be independent and exogenous. If we let $J_h^i:\mathbb{I}_T^i \times \mathbb{Y}^{[0,T)} \to \mathbb{R}$ denote DM$^i$'s sample-path cost function that determines DM$^i$'s total discounted cost realized at DM$^i$'s history $(I_T^i,Y_{T-1})$, then we have, for all $i\in[1,N]$, $\tilde{\gamma}^i\in\Gamma^i$, $\gamma\in\Gamma$,
$$E^{(\tilde{\gamma}^i,\zeta_{\gamma})}[J^i] = E^{\mu^{i,(\tilde{\gamma}^i,\zeta_{\gamma})}}[J_h^i] \quad\textrm{and}\quad
E^{(\tilde{\gamma}^i,\gamma^{-i})}[\bar{J}^{i,N}] = E^{\mu^{i,(\tilde{\gamma}^i,\gamma^{-i})}}[J_h^i].$$
We note that $J_h^i$ is independent of $N$, i.e., it is the same function in all games $\mathbb{G}^i,\mathbb{G}^{i+1},\dots$.
To show (\ref{eq:anonnashapprox}), we first show that, as $N\rightarrow\infty$,  the distance between $\mu^{i,(\tilde{\gamma}^i,\zeta_{\gamma})}$ and $\mu^{i,(\tilde{\gamma}^i,\gamma^{-i})}$, with respect to the Prokhorov metric\footnote{For a metric space $\mathbb{M}$, the Prokhorov metric on $\mathcal{P}(\mathbb{M})$ is defined by
$$m_{\mathcal{P}(\mathbb{M})}(\mu,\nu):=\inf\{\epsilon>0:\mu(B) \leq \nu(B_{\epsilon})+\epsilon, \ \forall B\in\mathcal{B}(\mathbb{M})\}, \qquad \forall \mu,\nu\in\mathcal{P}(\mathbb{M})$$ where $B_{\epsilon}:=\{y\in\mathbb{M}: m_{\mathbb{M}}(x,y)<\epsilon, \ \textrm{for some} \ x\in B \}$ denotes
the $\epsilon$ neighborhood of $B$.} that generates the weak topology, approaches zero uniformly for all $i\in[1,N]$, $\tilde{\gamma}^i\in\Gamma^i$, $\gamma\in\Gamma$.

For $i\in[1,N]$, $s=(s^i,s^{-i})\in\mathbb{S}$, define the mappings $F_s^0$ and $F_{s^i}^i$ by
$$Y_{T-1}=F_s^0(\omega) \quad \textrm{and} \qquad I_T^i=F_{s^i}^i(Y_{T-1},\omega^i)$$
where $\omega=(\omega^0,\omega^1,\dots,\omega^N)\in\Omega:=\times_{j\in[0,N]}\Omega^j$.
In other words, $F_s^0$ is the mapping that generates the environment variables $Y_{T-1}$ under the joint pure strategy $s$ and all of the random parameters $\omega$ in the system; whereas $F_{s^i}^i$ is the mapping that generates DM$^i$'s state-control trajectory $I_T^i$ under DM$^i$'s pure strategy $s^i$, DM$^i$'s random parameters $\omega^i$, and the given environment variables $Y_{T-1}$. $F_s^0$ is dependent on $N$, however, $F_{s^i}^i$ is not.
Let $\tilde{s}^i\in\mathbb{S}^i$ be drawn independently from $\tilde{\gamma}^i\in\Gamma^i$ ($s\in\mathbb{S}$ is drawn from $\Gamma$), and let $\tilde{\omega}^i$ be drawn independently from the distribution of $\omega^i$ (recall that $\omega^0,\omega^1,\dots,\omega^N$ are assumed to be mutually independent).
For any $B^i\in\mathcal{B}(\mathbb{I}_T^i \times \mathbb{Y}^{[0,T)})$, we have
\begin{align*}
\mu^{i,(\tilde{\gamma}^i,\zeta_{\gamma})}(B^i)= & P^{(\tilde{\gamma}^i,\gamma)}\big[(F_{\tilde{s}^i}^i(F_s^0(\omega),\tilde{\omega}^i),F_s^0(\omega))\in B^i\big] \\
\mu^{i,(\tilde{\gamma}^i,\gamma^{-i})}(B^i)= & P^{(\tilde{\gamma}^i,\gamma)}\big[( F_{\tilde{s}^i}^i(F_{(\tilde{s}^i,s^{-i})}^0(\tilde{\omega}^i,\omega^{-i}),\tilde{\omega}^i),F_{(\tilde{s}^i,s^{-i})}^0(\tilde{\omega}^i,\omega^{-i}))\in B^i\big]
\end{align*}
where $\omega^{-i}=(\omega^j)_{j\in[0,N]\backslash\{i\}}$.

Due to conditions of Assumption~\ref{as:anon1}, we have: given $\epsilon>0$, there exist $N_{\epsilon}\in\mathbb{N}$ such that, for all $N\geq N_{\epsilon}$, $i\in[1,N]$, $(\tilde{s}^i,\tilde{\omega}^i,s,\omega)\in\mathbb{S}^i\times\Omega^i\times\mathbb{S}\times\Omega$,
\begin{align*}
& m_{\mathbb{I}_T^i \times \mathbb{Y}^{[0,T)}} \big( ( F_{\tilde{s}^i}^i(F_s^0(\omega),\tilde{\omega}^i) , F_s^0(\omega)) \ , \ ( F_{\tilde{s}^i}^i(F_{(\tilde{s}^i,s^{-i})}^0(\tilde{\omega}^i,\omega^{-i}),\tilde{\omega}^i)) , F_{(\tilde{s}^i,s^{-i})}^0(\tilde{\omega}^i,\omega^{-i}) \big) \\ & \quad < \epsilon.
\end{align*}
As a result, given $\epsilon>0$, we have: for all $N\geq N_{\epsilon}$, $i\in[1,N]$, $(\tilde{s}^i,\tilde{\omega}^i,s,\omega)\in\mathbb{S}^i\times\Omega^i\times\mathbb{S}\times\Omega$, $B^i\in\mathcal{B}(\mathbb{I}_T^i\times\mathbb{Y}^{[0,T)})$,
\begin{align*}
( F_{\tilde{s}^i}^i(F_s^0(\omega),\tilde{\omega}^i) , F_s^0(\omega) )\in B^i & \  \Rightarrow \
( F_{\tilde{s}^i}^i(F_{(\tilde{s}^i,s^{-i})}^0(\tilde{\omega}^i,\omega^{-i}),\tilde{\omega}^i) , F_{(\tilde{s}^i,s^{-i})}^0(\tilde{\omega}^i,\omega^{-i}) )\in B_{\epsilon}^i
%\\
%\big\{( F_{\tilde{s}^i}^i(F_s^0(\omega),\tilde{\omega}^i) , F_s^0(\omega) )\in B_{\epsilon}^i\big\} & \supset
%\big\{( F_{\tilde{s}^i}^i(F_{(\tilde{s}^i,s^{-i})}^0(\tilde{\omega}^i,\omega^{-i}),\tilde{\omega}^i) , F_{(\tilde{s}^i,s^{-i})}^0(\tilde{\omega}^i,\omega^{-i}) )\in B^i\big\}
\end{align*}
therefore
\begin{align*}
\mu^{i,(\tilde{\gamma}^i,\zeta_{\gamma})}(B^i) \leq  \mu^{i,(\tilde{\gamma}^i,\gamma^{-i})}(B_{\epsilon}^i)
\end{align*}
where $B_{\epsilon}^i$ is the $\epsilon$ neighborhood of $B^i$.
This implies that
\begin{equation}
\label{eq:pmc0}
\lim_{N\rightarrow\infty} \sup_{i\in[1,N], \tilde{\gamma}^i\in\Gamma^i, \gamma\in\Gamma} m_{\mathcal{P}(\mathbb{I}_T^i \times \mathbb{Y}^{[0,T)})} (  \mu^{i,(\tilde{\gamma}^i,\zeta_{\gamma})} , \mu^{i,(\tilde{\gamma}^i,\gamma^{-i})}     ) =0.
\end{equation}
where $m_{\mathcal{P}(\mathbb{I}_T^i \times \mathbb{Y}^{[0,T)})}$ is the Prokhorov metric that generates the weak topology.

From (\ref{eq:pmc0}), we obtain, for each $j\in\mathbb{N}$,
\begin{equation}
\label{eq:wc0}
\lim_{N\rightarrow\infty} \sup_{i\in[1,N],\tilde{\gamma}^i\in\Gamma^i,\gamma\in\Gamma} |E^{\mu^{i,(\tilde{\gamma}^i,\zeta_{\gamma})}} [ J_h^j ] - E^{\mu^{i,(\tilde{\gamma}^i,\gamma^{-i})}} [ J_h^j ] | =0.
\end{equation}
We now observe that
$(J_h^i)_{i\in\mathbb{N}}$ is uniformly bounded and equi-continuous due to Condition~1(a)-(c) in Assumption~\ref{as:anon1}.
By Arzela-Ascoli Theorem \cite{Dud02}, $(J_h^i)_{i\in\mathbb{N}}$ is totally bounded for the uniform norm; therefore, it can be approximated by some $(J_h^i)_{i\in I}$, where $I$ is a finite subset of $\mathbb{N}$, with arbitrary precision with respect to the uniform norm.
As a result, (\ref{eq:wc0}) implies (\ref{eq:anonnashapprox}).

\subsection{Large  games with countably infinite set of DMs}
We introduce a large game with a countably infinite set of DMs, denoted by $\mathbb{G}^{\infty}$, based on a sequence $(\mathbb{G}^N)_{N\in\mathbb{N}}$ of games.
We consider finite as well as infinite time horizons, i.e., $T\in\mathbb{N}\cup\{\infty\}$.
Let $\mathbb{S}:=\times_{i\in\mathbb{N}} \mathbb{S}^i$ and $\Gamma:=\times_{i\in\mathbb{N}} \Gamma^i$.
The long-term cost of DM$^i$, $i\in\mathbb{N}$, in $\mathbb{G}^{\infty}$ corresponding to a joint mixed strategy $\gamma=(\gamma^j)_{j\in\mathbb{N}}\in\Gamma$ is defined by
\begin{equation}
\label{eq:barJinfdef}
\bar{J}^{i,\infty}(\gamma) = \limsup_{N\rightarrow\infty} E^{(\gamma^1,\dots,\gamma^N)} [\bar{J}^{i,N}]
\end{equation}
where $E^{(\gamma^1,\dots,\gamma^N)} [\bar{J}^{i,N}]$ denotes the long-term cost of DM$^i$ in $\mathbb{G}^N$, $N\in[i,\infty)$, corresponding to the joint mixed strategy $(\gamma^j)_{j\in[1,N]}$. A joint strategy $(\bar{\gamma}^j)_{j\in\mathbb{N}}\in  \Gamma$ satisfying $\bar{J}^{i,\infty}(\bar{\gamma}^i,\bar{\gamma}^{-i}) \leq \bar{J}^{i,\infty}(\gamma^i,\bar{\gamma}^{-i})$, for all $i\in\mathbb{N}$, $\gamma^i\in\Gamma^i$, constitutes a Nash equilibrium in mixed strategies in $\mathbb{G}^{\infty}$. The goal of this subsection is to show that some of the Nash equilibria in $\mathbb{G}^{\infty}$ can be obtained as the limiting SEBEU in $(\mathbb{G}^N)_{N\in\mathbb{N}}$ as $N\rightarrow\infty$.

We note that our results can be viewed as SEBEU counterparts of some related results in the literature which focus either on Nash equilibria or team optimality, but which also either impose or show the optimality of symmetric policies (unlike the setup we study here):
\cite{lasry2007mean,fischer2017connection,lacker2018convergence,bardi2014linear,arapostathis2017solutions,sanjari2018optimal,sanjari2019optimal,sanjari2020optimality}.

%By an abuse of notation, let us view any joint (pure or mixed) strategy in $\mathbb{G}^N$, $N\in\mathbb{N}$, as an element of $\mathbb{S}$ or $\Gamma$, whichever is appropriate, by appending to it some arbitrary element of $\times_{i> N}\mathbb{S}^i$ or $\times_{i> N}\Gamma^i$.
%Moreover, let us view any DM$^i$'s cost function $\bar{J}^i$ in $\mathbb{G}^N$, $N\in\mathbb{N}$, as a mapping from $\mathbb{S}$ to $\mathbb{R}$, by dropping the strategies for DM$^{N+1}$, DM$^{N+2}$,....

\begin{theorem}
\label{th:exanon1inf}
Consider a sequence $(\mathbb{G}^N)_{N\in\mathbb{N}}$ of games under Assumption~\ref{as:anon1}.
Let $(\gamma^{*,N})_{N\in\mathbb{N}}$ be such that $\gamma^{*,N}=(\gamma^{1,N},\dots,\gamma^{N,N})$ is an SEBEU in mixed strategies in $\mathbb{G}^N$ for each $N\in\mathbb{N}$. Suppose that there exists a joint mixed strategy $\gamma^{*,\infty}=(\gamma^{1,\infty},\gamma^{2,\infty},\dots)$ in $\mathbb{G}^{\infty}$ such that
\begin{equation}
\label{eq:pmc1}
\lim_{N\rightarrow\infty}  m_{\mathcal{P}(\mathbb{I}_T^i \times \mathbb{Y}^{[0,T)})} \big(\mu^{i,\gamma^{*,N}} , \mu^{i,(\gamma^{1,\infty},\dots,\gamma^{N,\infty})} \big) = 0, \quad \forall i\in\mathbb{N}
\end{equation}
where
$\mu^{i,(\gamma^{1,\infty},\dots,\gamma^{N,\infty})}$ denotes the probability measure induced on  $\mathcal{B}(\mathbb{I}_T^i \times \mathbb{Y}^{[0,T)})$ by the joint mixed strategy $(\gamma^{1,\infty},\dots,\gamma^{N,\infty})$ in $\mathbb{G}^N$.
Then, $\gamma^{*,\infty}$ is a Nash equilibrium in $\mathbb{G}^{\infty}$.
\end{theorem}

\textbf{Proof of Theorem~\ref{th:exanon1inf}}
By Theorem~\ref{th:exanon1}, given $\epsilon>0$, there exist $\bar{N}_{\epsilon}\in\mathbb{N}$ such that, if $N\geq \bar{N}_{\epsilon}$, then
$\gamma^{*,N}$ is an $\epsilon$-Nash equilibrium in $\mathbb{G}^N$, i.e.,
\begin{equation}
\label{eq:pmc2}
E^{\gamma^{*,N}} [\bar{J}^{i,N}] \leq E^{(\tilde{\gamma}^{i},\gamma^{-i,N})} [\bar{J}^{i,N}] + \epsilon, \quad \forall i\in[1,N],  \tilde{\gamma}^i\in\Gamma^i.
\end{equation}
By (\ref{eq:pmc1}), we have
$$
\lim_{N\rightarrow\infty} |E^{\gamma^{*,N}}[\bar{J}^{i,N}] - E^{(\gamma^{1,\infty},\dots,\gamma^{N,\infty})}[\bar{J}^{i,N}]| = 0, \quad \forall i\in\mathbb{N}
$$
which, together with (\ref{eq:barJinfdef}), implies that
\begin{equation}
\label{eq:pmc3}
\bar{J}^{i,\infty}(\gamma^{*,\infty})
= \limsup_{N\rightarrow\infty} E^{\gamma^{*,N}} [\bar{J}^{i,N}], \quad \forall i\in\mathbb{N}.
\end{equation}
By (\ref{eq:pmc1}), we also have
$$\lim_{N\rightarrow\infty}  |E^{(\tilde{\gamma}^i,\zeta_{\gamma^{*,N}})}[J^i] - E^{(\tilde{\gamma}^i,\zeta_{(\gamma^{1,\infty},\dots,\gamma^{N,\infty})})}[J^i]| = 0, \quad \forall i\in\mathbb{N},  \tilde{\gamma}^i\in\Gamma^i.$$
In view of (\ref{eq:anonnashapprox}), this implies that
$$\lim_{N\rightarrow\infty}  |E^{(\tilde{\gamma}^i,\gamma^{-i,N})}[\bar{J}^{i,N}] - E^{(\tilde{\gamma}^i,(\gamma^{j,\infty})_{j\in[1,N]\backslash\{i\}})}[\bar{J}^{i,N}]| = 0, \quad \forall i\in\mathbb{N},  \tilde{\gamma}^i\in\Gamma^i$$
which, together with (\ref{eq:barJinfdef}), implies
\begin{equation}
\label{eq:pmc4}
\bar{J}^{i,\infty}(\tilde{\gamma}^i,\gamma^{-i,\infty})
= \limsup_{N\rightarrow\infty} E^{(\tilde{\gamma}_i,\gamma^{-i,N})} [\bar{J}^{i,N}], \quad \forall i\in\mathbb{N},  \tilde{\gamma}^i\in\Gamma^i.
\end{equation}
From (\ref{eq:pmc2})-(\ref{eq:pmc4}), we obtain, for any $\epsilon>0$,
$$\bar{J}^{i,\infty}(\gamma^{*,\infty})  \leq \bar{J}^{i,\infty}(\tilde{\gamma}^i,\gamma^{-i,\infty}) + \epsilon, \quad \forall i\in[1,N],  \tilde{\gamma}^i\in\Gamma^i.
$$
Since $\epsilon>0$ is arbitrary, the desired result follows.

\subsection{SEBEU in pure strategies}
Under additional assumptions, it is possible to show the existence of $\epsilon-$SEBEU in pure strategies in large games. These pure strategies, however, are typically different for different DMs, i.e., they do not constitute a symmetric $\epsilon-$SEBEU.
\begin{assumption}
\label{as:anon2}
$\textrm{}$

\begin{enumerate}
\item DMs are identical with $\mathbb{X}^1  \subset \mathbb{R}^n$ and $\mathbb{U}^1  \subset \mathbb{R}^m$ for some $n,m\in\mathbb{N}$
\item The aggregate statistics are the simple averages in (\ref{eq:avgfun})
\item The environment variables $(y_t)_{t\in[0,T)}$ are generated by deterministic dynamics driven by the (random) inputs $(d_t^N)_{t\in[0,T)}$, i.e.,
$\Omega^0$ is a singleton and, by a slight abuse of notation,
\begin{align}
\label{eq:p3}
y_t =  & g_t(x_t^0, d_t^N), \qquad t\in[0,T) \\
\label{eq:y3}
x_{t+1}^0 = & f_t^0(x_t^0, d_t^N), \qquad t\in[0,T).
\end{align}
\end{enumerate}
\end{assumption}

\begin{theorem}
\label{th:identpure1}
Consider a sequence $(\mathbb{G}^N)_{N\in\mathbb{N}}$ of games under Assumption~\ref{as:anon1}-\ref{as:anon2}.
Given $\epsilon>0$, there exists $\hat{N}_{\epsilon}\in\mathbb{N}$ such that, if $N\geq \hat{N}_{\epsilon}$, then $\mathbb{G}^N$ possess an $\epsilon-$SEBEU in pure strategies, i.e., there exists $s=(s^1,\dots,s^N)\in\mathbb{S}$ such that, for all $i\in[1,N]$, $\tilde{\gamma}^i\in\Gamma^1$, $$E^{\zeta_s} [ J^i(s^i,\cdot) ] \leq E^{(\tilde{\gamma}^i,\zeta_s)} [ J^i ] +\epsilon.$$
\end{theorem}

Such pure-strategy $\epsilon-$SEBEU can be constructed through purification of symmetric mixed-strategy  SEBEU $\gamma^{*.N}$ in $\mathbb{G}^N$ (which exists by Theorem~\ref{thm:cm}) in the spirit of \cite{schmeidler1973equilibrium}. The proof of Theorem~\ref{th:identpure1}, presented in \ref{appx:pf}, relies on the following intuition. As the number of DMs increases, the aggregate statistics $(d_t^N)_{t\in[0,T)}$ and the environment variables $(y_t)_{t\in[0,T)}$ generated under $\gamma^{*,N}$ converge in probability to suitably defined ``mean values''. This is due to the mutual independence of the local random parameters $(\omega^i)_{i\in\mathbb{N}}$, the identicalness of the DMs, the deterministic nature of the dynamics (\ref{eq:p3})-(\ref{eq:y3}) generating the environment variables, and the continuity of the system functions (Assumption~\ref{as:anon1}). The key observation is that $\gamma^{*,N}$ can be approximated by a joint pure strategy $s^{*,N}$ in such a way that the aggregate statistics $(d_t^N)_{t\in[0,T)}$ and the environment variables $(y_t)_{t\in[0,T)}$ generated under $s^{*,N}$ also converge to the same mean values as the number of DMs increases.
This allows us to show that $s^{*,N}$ constitutes a pure-strategy $\epsilon-$SEBEU in $\mathbb{G}^N$ when the number of DMs is sufficiently large.
%The details of the proof of Theorem~\ref{th:identpure1} is provided in Appendix~\ref{appx:anon}.

The assumption of identical DMs in Theorem~\ref{th:identpure1} can be relaxed to different sets of identical DMs (DMs in different sets need not be identical), e.g., when the number of identical DMs within each set in $\mathbb{G}^N$ grows  without bound as $N\rightarrow\infty$.
%The fourth condition of Assumption~\ref{as:anon1}, namely the independence of the local random variables $\omega^1,\omega^2,\dots$, can also be relaxed as long as the weak law of large numbers apply to the samples $((x_t^i,u_t^i))_{i\in\mathcal{I}^N(s^1)}$, for each $s^1\in\mathbb{S}^1$ and $t\in[0,T)$, as $N\rightarrow\infty$.
We should also point out that, for sufficiently large number of DMs, Theorem~\ref{th:identpure1} together with Theorem~\ref{th:exanon1} imply the existence of pure-strategy $\epsilon-$Nash equilibrium  in games satisfying the conditions of Theorem~\ref{th:identpure1}.

\section{Conclusion}
We presented a generalization of the concept of price-taking equilibrium to general stochastic dynamic games, called SEBEU.
At an SEBEU, players optimize their long-term cost  with respect to their beliefs of exogenous uncertainty in the environment variables while their beliefs being consistent with the true distribution of the environment variables.
We established general existence of SEBEU in mixed and behavior strategies when all state, control, disturbance, and environment variables belong to compact metric spaces.
We showed the near person-by-person optimality of SEBEU in large games with small players. Devising learning methods that lead to SEBEU in stochastic dynamic games remains as a significant future research topic.

% Appendix here
% Options are (1) APPENDIX (with or without general title) or
%             (2) APPENDICES (if it has more than one unrelated sections)
% Outcomment the appropriate case if necessary
%

\appendix

\section{Proof of existence of SEBEU in mixed strategies for compact metric space models}
\label{appx:ex}
We provide a proof of Theorem~\ref{thm:cm} by verifying that the mapping $\overline{BR}:\Gamma\rightrightarrows\Gamma$ satisfies the conditions of Fan-Glicksberg's fixed point theorem \cite{Fan1952,Glicksberg1952} and \cite{fud-tirole}:
\begin{enumerate}
\item $\Gamma$ is a compact, convex, nonempty subset of a locally convex topological vector space.
\item $\overline{BR}(\gamma)$ is nonempty for all $\gamma$.
\item $\overline{BR}(\gamma)$ is convex for all $\gamma$.
\item $\overline{BR}(\cdot)$ has a closed graph, i.e., if $(\gamma_n,\check{\gamma}_n)\rightarrow(\gamma,\check{\gamma})$ with $\check{\gamma}_n \in \overline{BR} (\gamma_n)$, then $\check{\gamma} \in \overline{BR} (\gamma)$.
\end{enumerate}

We now verify each of these conditions.
The set $\mathbb{S}^i$ of pure strategies for any DM$^i$ is a countable product of the compact metric spaces $(\mathbb{S}_t^i)_{t\in[0,T)}$, hence $\mathbb{S}^i$ a compact\footnote{The cartesian product of any nonempty collection of compact spaces is compact in the product topology due to the Tychnonoff Theorem.} space metrizable by the product metric in (\ref{eq:pmet}). This implies that the set $\Gamma^i$ of mixed strategies is a nonempty, convex, compact subset of a locally convex topological vector space, namely the space of finite signed Borel measures on $\mathbb{S}^i$ endowed with the weak topology\footnote{The weak topology is the coarsest topology under which the mappings $\gamma^i \mapsto E^{\gamma^i} [f]$ are continuous for any bounded continuous function $f:\mathbb{S}^i\to\mathbb{R}$.}; see Theorem~6.4 in \cite{parthasarathy2014probability}.
Hence, the first condition holds.

The mapping $(s^i,Z_{T-1})\mapsto J^i(s^i,Z_{T-1})$ is continuous in $\mathbb{S}^i \times \mathbb{Y}^{[0,T)}$ under the product topology.
When $T\in\mathbb{N}$, this follows from the uniform continuity of the functions $c_0^i,\dots,c_T^i,f_0^i,\dots,f_{T-1}^i$ and the uniform equi-continuity of the mappings in $\mathbb{S}_0^i,\dots,\mathbb{S}_{T-1}^i$.
When $T=\infty$, due to Assumption~\ref{as:cm}, for any $\epsilon>0$, there exists $T_{\epsilon/4}\in\mathbb{N}$ such that, the long-term cost incurred over the horizon $t\in[T_{\epsilon/4},\infty)$, under any $(s^i,Z_{T-1})\in\mathbb{S}^i\times \mathbb{Y}^{[0,T)}$, satisfies
\begin{equation}
\label{eq:costtrunc}\Bigg|E \Bigg[\sum_{t\in[T_{\epsilon/2},\infty)} (\beta^i)^t c_t^i\big(x_t^i,u_t^i,z_t\big)   \Bigg]\Bigg|< \epsilon/4
\end{equation}
and in addition there exists $\delta>0$ such that
$$
m_{\mathbb{S}^i\times \mathbb{Y}^{[0,T)}}((s^i,Z_{T-1}),(\bar{s}^i,\bar{Z}_{T-1})) < \delta \quad\Rightarrow\quad
|J_{T_{\epsilon/4}}^i(s^i,Z_{T-1})-J_{T_{\epsilon/4}}^i(\bar{s}^i,\bar{Z}_{T-1})|<\epsilon/2
$$
where $J_{T_{\epsilon/4}}^i(s^i,Z_{T-1})$ represents DM$^i$'s (truncated) cost incurred over the horizon $t\in[0,T_{\epsilon/4})$ corresponding to $(s^i,Z_{T-1})$.
This implies the continuity of the mapping $(s^i,Z_{T-1})\mapsto J^i(s^i,Z_{T-1})$ in $\mathbb{S}^i\times \mathbb{Y}^{[0,T)}$. Since the mapping $(s^i,Z_{T-1})\mapsto J^i(s^i,Z_{T-1})$ is also bounded in $\mathbb{S}^i\times \mathbb{Y}^{[0,T)}$ due to Assumption~\ref{as:cm}, the mapping $(\gamma^i,\zeta)\mapsto E^{(\gamma^i,\zeta)} [J^i]$ is continuous in $\Gamma^i\times\mathcal{P}(\mathbb{Y}^{[0,T)})$ under the weak topology.
The second condition now follows from the fact that each $\overline{BR}^i(\gamma)$ is the set of minimizers of the continuous function $E^{(\cdot,\zeta_{\gamma})} [J^i]$ over the nonempty compact set $\Gamma^i$.

For any $\zeta\in\mathcal{P}(\mathbb{Y}^{[0,T)})$, the mapping $\gamma^i\mapsto E^{(\gamma^i,\zeta)} [J^i]$ is clearly convex, which implies the third condition.

To verify the fourth condition, assume that $(\gamma_n,\check{\gamma}_n)\rightarrow(\gamma,\check{\gamma})$ with $\check{\gamma}_n \in \overline{BR} (\gamma_n)$, but $\check{\gamma} \not\in \overline{BR} (\gamma)$. This means that, for some $i\in[1,N]$, $\tilde{\gamma}^i\in\Gamma^i$, $\epsilon>0$, we have $E^{(\tilde{\gamma}^i,\zeta_{\gamma})} [J^i] < E^{(\check{\gamma}^i,\zeta_{\gamma})} [J^i] - 3\epsilon$.
We note that the mapping $(\tilde{s}^i,s) \mapsto E^{(\tilde{s}^i,\zeta_s)} [J^i]$, where $E^{(\tilde{s}^i,\zeta_s)} [J^i]=E^{\zeta_s} [J^i(\tilde{s}^i,\cdot)]$, is continuous and bounded in $\mathbb{S}^i\times\mathbb{S}$ due to  Assumption~\ref{as:cm}, in particular, due to (\ref{eq:costtrunc}) and the uniform equi-continuity of the mappings in each $\mathbb{S}^i$. This implies the continuity of the mapping $(\tilde{\gamma}^i,\gamma) \mapsto E^{(\tilde{\gamma}^i,\zeta_{\gamma})} [J^i]$ in $\Gamma^i\times\Gamma$ because, for any $\epsilon>0$,
$$|E^{(\tilde{\gamma}^i,\gamma)} [E^{(\tilde{s}^i,\zeta_s)} [J^i]] - E^{(\hat{\gamma}^i,\bar{\gamma})} [E^{(\tilde{s}^i,\zeta_s)} [J^i]]|<\epsilon \quad\Rightarrow\quad
|E^{(\tilde{\gamma}^i,\zeta_{\gamma})} [J^i] - E^{(\hat{\gamma}^i,\zeta_{\bar{\gamma}})} [J^i]| < \epsilon$$
where the inner expectation in $E^{(\tilde{\gamma}^i,\gamma)} [ E^{(\tilde{s}^i,\zeta_s)} [ J^i ] ]$ is with respect to the distribution $\zeta_s$ induced by a fixed pure strategy profile $s\in\mathbb{S}$ (recall $E^{(\tilde{s}^i,\zeta_s)} [J^i]=E^{\zeta_s} [J^i(\tilde{s}^i,\cdot)$)] and the outer expectation is with respect to the probability distribution $(\tilde{\gamma}^i,\gamma)$ over the set $\mathbb{S}^i\times\mathbb{S}$ of pure strategies.
Therefore, for sufficiently large $n\in\mathbb{N}$, we have
$$  E^{(\tilde{\gamma}^i,\zeta_{\gamma_n})} [J^i] -\epsilon  <  E^{(\tilde{\gamma}^i,\zeta_{\gamma})} [J^i] < E^{(\check{\gamma}^i,\zeta_{\gamma})} [J^i] - 3\epsilon < E^{(\check{\gamma}_n^i,\zeta_{\gamma_n})} [J^i] - 2\epsilon.$$
This contradicts $\check{\gamma}_n^i \in \overline{BR}^i(\gamma_n)$; hence, the fourth condition also holds. The existence of an SEBEU in mixed strategies is now implied by Fan-Glicksberg's fixed point theorem.

If there are sets of identical DMs,
the existence of a symmetric SEBEU in mixed strategies can be shown similarly by applying Fan-Glicksberg's fixed theorem on the correspondence $\gamma \mapsto (BR^1(\zeta_{\gamma}),\dots,BR^N(\zeta_{\gamma}))$ where $\gamma\in\Gamma$ is restricted in such a way that identical DMs have identical strategies at $\gamma$.

\section{Proof of Theorem~\ref{th:identpure1}}
\label{appx:pf}
%\textbf{Proof of Theorem~\ref{th:identpure1}}
Due to (\ref{eq:costtruncinf}), it suffices to prove the theorem for the case $T\in\mathbb{N}$. By Theorem~\ref{thm:cm}, $\mathbb{G}^N$ possesses a symmetric SEBEU $\gamma^{*,N}=(\gamma^{1,N},\cdots,\gamma^{N,N}) \in\Gamma$, where $\gamma^{1,N}=\cdots=\gamma^{N,N}$.
%Let $d_t^N$ and $y_t^N$ denote, respectively, the aggregate statistic and the environment variable at time $t\in[0,T)$ generated in $\mathbb{G}^N$ under some joint strategy.
%Let $\bar{d}_t^N$ and $\bar{y}_t^N$, $t\in[0,T)$, be generated  by the following ``mean'' dynamics, for all $i\in[1,N]$, $t\in[0,T)$,
%\begin{align}
%\bar{d}_t^N = & \frac{1}{N}\sum_{i\in[1,N]} \left(\begin{array}{c} \bar{x}_t^i \\ \bar{u}_t^i \end{array} \right) \nonumber \\
%\bar{y}_t^N = & g_t(\bar{x}_t^0,E^{\gamma^{*,N}} [\bar{d}_t^N] ) \label{eq:barytN}\\
%\bar{x}_{t+1}^0 = & f_t^0(\bar{x}_t^0,E^{\gamma^{*,N}} [\bar{d}_t^N]) \label{eq:barx0tN}\\
%\bar{x}_{t+1}^i = & f_t^i(\bar{x}_t^i,\bar{u}_t^i,\bar{y}_t^N,w_t^i) \nonumber  \\
%\bar{u}_t^i     = & s_t^i(\bar{I}_t^i,\bar{Y}_{t-1}) \nonumber
%\end{align}
Let $(y_t^N)_{t\in[0,T)}$, $(\bar{y}_t^N)_{t\in[0,T)}$ be generated in $\mathbb{G}^N$ by the following dynamics: for all $i\in[1,N]$, $t\in[0,T)$,
\begin{align}
d_t^N = & \frac{1}{N}\sum_{i\in[1,N]} \begin{bmatrix} x_t^i \\ u_t^i \end{bmatrix},
& \bar{d}_t^N = & \frac{1}{N}\sum_{i\in[1,N]} \begin{bmatrix} \bar{x}_t^i \\ \bar{u}_t^i \end{bmatrix}\label{eq:bardtN} \\
y_t^N = & g_t(x_t^0,d_t^N ),
& \bar{y}_t^N = & g_t(\bar{x}_t^0,E^{\gamma^{*,N}} [\bar{d}_t^N] )\label{eq:barytN}\\
x_{t+1}^0 = & f_t^0(x_t^0,d_t^N),
& \bar{x}_{t+1}^0 = & f_t^0(\bar{x}_t^0,E^{\gamma^{*,N}} [\bar{d}_t^N])\label{eq:barx0tN}\\
x_{t+1}^i = & f_t^i(x_t^i,u_t^i,y_t^N,w_t^i),
& \bar{x}_{t+1}^i = & f_t^i(\bar{x}_t^i,\bar{u}_t^i,\bar{y}_t^N,w_t^i) \label{eq:barxtN}  \\
u_t^i     = & s_t^i(I_t^i,Y_{t-1}),
& \bar{u}_t^i     = & s_t^i(\bar{I}_t^i,\bar{Y}_{t-1})\label{eq:barstN}
\end{align}
where $(\bar{x}_0^0,\bar{x}_0^1,\dots,\bar{x}_0^N)=(x_0^0,x_0^1,\dots,x_0^N)$, $s^i=(s_t^i)_{t\in[0,T)}\in\mathbb{S}^1$ is some pure strategy for DM$^i$, and $(\bar{I}_t^i,\bar{Y}_{t-1})$ is defined analogously to $(I_t^i,Y_{t-1})$ (we suppressed the dependence of the various terms on $N$ other than $d_t^N$, $y_t^N$, $\bar{d}_t^N$, $\bar{y}_t^N$).
We have a few observations about the dynamics (\ref{eq:bardtN})-(\ref{eq:barstN}).
\begin{itemize}
\item[-] The same pure strategies $s=(s^i)_{i\in[1,N]}$ and the same realizations of the random parameters $\omega\in\Omega$ generate
$\big((\bar{d}_t^N,\bar{y}_t^N,\bar{x}_t^0,(\bar{x}_t^i)_{i\in[1,N]},(\bar{u}_t^i)_{i\in[1,N]})\big)_{t\in[0,T)}$ and $\big((d_t^N,y_t^N,x_t^0,(x_t^i)_{i\in[1,N]},(u_t^i)_{i\in[1,N]})\big)_{t\in[0,T)}$.
\item[-] The joint pure strategy $s=(s^i)_{i\in[1,N]}$ is possibly selected according to some joint mixed strategy, not necessarily $\gamma^{*,N}$.
\item[-] $(y_t^N)_{t\in[0,T)}$ represent the environment variables generated in $\mathbb{G}^N$.
\item[-] $(\bar{y}_t^N)_{t\in[0,T)}$ is a deterministic sequence generated by (\ref{eq:barytN})-(\ref{eq:barx0tN}), $\bar{x}_0^0$, and the inputs $(E^{\gamma^{*,N}}[\bar{d}_t^N])_{t\in[0,T)}$.
\item[-]  $(\bar{y}_t^N)_{t\in[0,T)}$ is not influenced by the joint pure strategy $s=(s^i)_{i\in[1,N]}$ in (\ref{eq:barstN}).
\end{itemize}
The proof consists of showing the following steps.
\begin{itemize}
\item[(i)] If the joint strategy $s=(s^i)_{i\in[1,N]}$ in (\ref{eq:barstN}) is selected according to $\gamma^{*,N}$, $|\bar{d}_t^N - E^{\gamma^{*,N}}[\bar{d}_t^N]|$ converges in probability to zero as $N\rightarrow0$, i.e.,
\begin{equation}
\label{eq:notcipdg}
\lim_{N\rightarrow\infty} P^{\gamma^{*,N}}[|\bar{d}_t^N - E^{\gamma^{*,N}}[\bar{d}_t^N]| > \delta]=0, \quad\forall \delta>0, t\in[0,T).
\end{equation}
\item[(ii)] There exist  $s^{*,N}\in\mathbb{S}$ such that
    \begin{equation}
\label{eq:notcipds}
\lim_{N\rightarrow\infty} P^{s^{*,N}}[|\bar{d}_t^N - E^{\gamma^{*,N}}[\bar{d}_t^N]| > \delta]=0, \quad\forall \delta>0, t\in[0,T).
\end{equation}
\item[(iii)]  If the joint strategy $s=(s^i)_{i\in[1,N]}$ in (\ref{eq:barstN}) is selected according to $\gamma^{*,N}$ or as $s^{*,N}$, $m_{\mathbb{Y}} (y_t^N,\bar{y}_t^N)$ converges in probability to zero as $N\rightarrow0$, i.e.,
\begin{align}
\label{eq:notcip}
\lim_{N\rightarrow\infty} P^{\gamma^{*,N}}[m_{\mathbb{Y}} (y_t^N,\bar{y}_t^N) > \delta] & = 0, \quad\forall\delta>0, t\in[0,T) \\
\lim_{N\rightarrow\infty} P^{s^{*,N}}[m_{\mathbb{Y}} (y_t^N,\bar{y}_t^N) > \delta] & = 0, \quad\forall\delta>0, t\in[0,T) .
\label{eq:notcip2}
\end{align}
\item[(iv)] $s^{*,N}$ is an approximate Nash equilibrium in $\mathbb{G}^N$ for sufficiently large $N\in\mathbb{N}$.
\end{itemize}
Below, we use the notation $\xrightarrow[N\rightarrow\infty]{\gamma^{*,N}}$ or $\xrightarrow[N\rightarrow\infty]{s^{*,N}}$ for convergence in probability as $N\rightarrow\infty$ under the joint strategies $(\gamma^{*,N})_{N\in\mathbb{N}}$ or $(s^{*,N})_{N\in\mathbb{N}}$, respectively, e.g., as in (\ref{eq:notcipdg})-(\ref{eq:notcip2}).
\\

Step (i): We assume in this step that $s=(s^i)_{i\in[1,N]}$ in (\ref{eq:barstN}) is selected according to $\gamma^{*,N}$.
We note that $(\bar{y}_t^N)_{t\in[0,T)}$ is a deterministic sequence and $\gamma^{1,N}=\cdots=\gamma^{N,N}$; therefore,
$\left(\left[\begin{smallmatrix} \bar{x}_t^i \\ \bar{u}_t^i \end{smallmatrix}\right]\right)_{i\in[1,N]}$ are mutually independent and identically distributed, for each $t\in[0,T)$.
Hence, from Chebyshev's inequality and the assumption that $\mathbb{X}^1$, $\mathbb{U}^1$ are compact metric spaces, we obtain (as desired)
\begin{equation}
P^{\gamma^{*,N}}[|\bar{d}_t^{N}-E^{\gamma^{*,N}} [\bar{d}_t^{N}]| >  \delta] \leq \ \frac{1}{N\delta^2}\mathrm{trace}\left(\mathrm{cov}^{\gamma^{1,N}}\begin{bmatrix}\bar{x}_t^1 \\ \bar{u}_t^1 \end{bmatrix}\right) \leq \frac{D^2}{N\delta^2}, \quad\forall \delta>0, t\in[0,T)
\label{eq:inequ1}
\end{equation}
where
$$D:=\sup_{\hat{x}^1,\check{x}_1\in\mathbb{X}^1,\hat{u}^1,\check{u}^1\in\mathbb{U}^1} \max\left(\left| \begin{bmatrix} \hat{x}^1 \\ \hat{u}^1 \end{bmatrix} - \begin{bmatrix} \check{x}^1 \\ \check{u}^1 \end{bmatrix} \right| , \left| \begin{bmatrix} \hat{x}^1 \\ \hat{u}^1 \end{bmatrix} \right| \right)<\infty.$$
Letting $N\rightarrow\infty$ in (\ref{eq:inequ1}) gives us (\ref{eq:notcipdg}).

Step (ii): We begin by constructing a joint pure strategy $s^{*,N}$  under which (\ref{eq:notcipds}) holds. Note that the common set $\mathbb{S}^1$ of pure strategies for every DM$^i$ in every game $\mathbb{G}^N$ is a compact metric space.
Hence, $\mathbb{S}^1$ can be covered by finitely many balls of arbitrarily small radius.
Since $f_0^1,\dots,f_{T-1}^1$ are uniformly continuous and the mappings in $\mathbb{S}^1$ are uniformly equi-continuous, for each $\rho>0$, there exist a minimum number $n(\rho)\in\mathbb{N}$ and non-empty disjoint Borel subsets  $B_{\rho}^1,\dots,B_{\rho}^{n(\rho)}\subset\mathbb{S}^1$ such that $\cup_{\ell\in[1,n(\rho)]} B_{\rho}^{\ell} = \mathbb{S}^1$ and
\begin{equation}
\label{eq:Dbound}
\sup_{\hat{s}^1,\check{s}^1\in B_{\rho}^{\ell}}  \left|E^{\hat{s}^1}\begin{bmatrix}\bar{x}_t^1 \\ \bar{u}_t^1\end{bmatrix} - E^{\check{s}^1} \begin{bmatrix}\bar{x}_t^1 \\ \bar{u}_t^1\end{bmatrix} \right| \leq \rho, \quad \forall \ell\in[1,n(\rho)], t\in[0,T)
\end{equation}
where $\left(\left[\begin{smallmatrix} \bar{x}_t^1 \\ \bar{u}_t^1 \end{smallmatrix}\right]\right)_{t\in[0,T)}$ is generated by (\ref{eq:bardtN})-(\ref{eq:barstN}).
Let us define
\begin{align*}
\rho^N & :=  \inf \{\rho > 0 : n(\rho) \leq 1 + \ln(N) \} + 1/N \quad\textrm{and}\quad
n^N     := n(\rho^N), \quad N\in\mathbb{N}.
\end{align*}
We have $$\lim_{N\rightarrow\infty} \rho^N=0 \qquad\textrm{and}\qquad \lim_{N\rightarrow\infty} \frac{n^N}{\sqrt{N}} =0.$$
For ease of notation, let us also introduce $B^{\ell,N}:=B_{\rho^N}^{\ell}$, for each $N\in\mathbb{N}$, $\ell\in[1,n^N]$.

Fix $N\in\mathbb{N}$.
For each $\ell\in[1,n^N]$ such that $\gamma^{1,N}(B^{\ell,N})>0$, let us fix a strategy $\hat{s}^{\ell} \in B^{\ell,N}$ satisfying
\begin{equation}
\label{eq:hatslopt}
E^{\zeta_{\gamma^{*,N}}} [ J^i(\hat{s}^{\ell},\cdot) ] = E^{(\gamma^{i,N},\zeta_{\gamma^{*,N}})} [ J^i ], \qquad \forall i\in[1,N].
\end{equation}
Such $\hat{s}^{\ell} \in B^{\ell,N}$ exists when $\gamma^{1,N}(B^{\ell,N})>0$ because $\gamma^{*,N}$ is an SEBEU.

We will introduce $s^{*,N}$ by assigning each DM to some $\hat{s}^{\ell}$, $\ell\in[1,n^N]$, such that the fraction of DMs
assigned to $\hat{s}^{\ell}$ is approximately equal to $\gamma^{1,N}(B^{\ell,N})$. To achieve this, let $k_0:=0$, and define the number $k_{\ell}$ of DMs assigned to $\hat{s}^{\ell}$, recursively,  as
\begin{align*}
k_{\ell} := & \left\lfloor N\sum_{j=1}^{\ell}\gamma^{1,N}(B^{j,N})+1/2\right\rfloor-\sum_{j=0}^{\ell-1} k_{j}, \qquad \ell\in[1,n^N]
\end{align*}
where $\lfloor\cdot\rfloor$ denotes the integer floor. We have $k_1,\dots,k_{n^N}\geq0$ and $\sum_{\ell=1}^{n^N} k_{\ell}=N$. Moreover,
\begin{align}
\gamma^{1,N}(B^{\ell,N})=0 \ \Rightarrow \ k_{\ell}=0   \qquad\textrm{and}\qquad \ \max_{\ell\in[1,n^N]}\bigg|\frac{k_{\ell}}{N}  - \gamma^{1,N}(B^{\ell,N})\bigg| \leq \frac{1}{N}. \label{eq:s*appx3}
\end{align}
We now define $s^{*,N}=(s^{1,N},\dots,s^{N,N})$ by assigning each of $(\textrm{DM}^i)_{i\in[1,k_1]}$ to $\hat{s}^1$, each of $(\textrm{DM}^i)_{i\in[k_1+1,k_1+k_2]}$ to $\hat{s}^2$, and so on, i.e.,
\begin{align}
s^{j,N} :=\hat{s}^{\ell}, \qquad j\in\mathcal{I}_{\ell}.
\label{eq:s*Nconstr}
\end{align}
where $\mathcal{I}_{\ell}:=[k_0+\cdots+k_{\ell-1}+1,k_0+\cdots k_{\ell}]$, $\ell\in[1,n^N]$.

We then assume that $s=(s^i)_{i\in[1,N]}$ in (\ref{eq:barstN}) is selected as $s^{*,N}$ (for the rest of Step (ii)) and show (\ref{eq:notcipds}).
We note that $\left(\left[\begin{smallmatrix}\bar{x}_t^i \\ \bar{u}_t^i\end{smallmatrix}\right]\right)_{i\in\mathcal{I}_{\ell}}$  are mutually independent and identically distributed, for each $\ell\in[1,n^N]$, $t\in[0,T)$.
Hence, from Chebyshev's inequality, we obtain, for each $t\in[0,T)$, $N\in\mathbb{N}$, $\ell\in[1,n^N]$  with $k_{\ell}>0$,
\begin{align}
& P^{s^{*,N}}\bigg[\bigg| \underbrace{\frac{1}{k_{\ell}} \sum_{i\in\mathcal{I}_{\ell}} \begin{bmatrix} \bar{x}_t^i \\ \bar{u}_t^i \end{bmatrix} - E^{\hat{s}^{\ell}} \begin{bmatrix} \bar{x}_t^1 \\ \bar{u}_t^1 \end{bmatrix}}_{:=\alpha_{t,\ell}} \bigg|  >  \delta \bigg]  \leq \frac{1}{k_{\ell}\delta^2} \mathrm{trace}\left(\mathrm{cov}^{\hat{s}^{\ell}}\begin{bmatrix}\bar{x}_t^1 \\ \bar{u}_t^1 \end{bmatrix}\right) \leq   \frac{D^2}{k_{\ell}\delta^2}.
\label{eq:alphabound}
\end{align}
We write $\bar{d}_t^N - E^{\gamma^{*,N}} [ \bar{d}_t^N ] $ as
\begin{align}
& \bar{d}_t^N - E^{\gamma^{*,N}} [ \bar{d}_t^N ]  \\ & \qquad =  \sum_{\ell\in[1,n^N]:k_{\ell}>0} \frac{k_{\ell}}{N} \alpha_{t,\ell} + \sum_{\ell\in[1,n^N]} \bigg(\frac{k_{\ell}}{N} - \gamma^{1,N}(B^{\ell,N}) \bigg) E^{\hat{s}^{\ell}} \begin{bmatrix}\bar{x}_t^1 \\ \bar{u}_t^1 \end{bmatrix}
\nonumber \\
& \qquad \quad + \sum_{\ell\in[1,n^N]:\gamma^{1,N}(B^{\ell,N})>0} \gamma^{1,N}(B^{\ell,N})  \left( E^{\hat{s}^{\ell}} \begin{bmatrix}\bar{x}_t^1 \\ \bar{u}_t^1 \end{bmatrix}- E^{\gamma^{1,N}} \left[\left.\begin{bmatrix}\bar{x}_t^1 \\ \bar{u}_t^1 \end{bmatrix} \right| B^{\ell,N}  \right] \right). \label{bardtNdiffs*}
\end{align}
We upper-bound the first sum on the right-hand-side of (\ref{bardtNdiffs*}) as
\begin{align*}
\sum_{\ell\in[1,n^N]:k_{\ell}>0} \frac{k_{\ell}}{N} |\alpha_{t,\ell}| & \leq \max_{\ell\in[1,n^N]:k_{\ell}>\sqrt{N}} |\alpha_{t,\ell}| + \sum_{\ell\in[1,n^N]:0<k_{\ell}\leq \sqrt{N}} \frac{k_{\ell}}{N} |\alpha_{t,\ell}| \\ & \leq
 \max_{\ell\in[1,n^N]:k_{\ell}>\sqrt{N}} |\alpha_{t,\ell}| + n^N D/\sqrt{N}
\end{align*}
where we set $\max_{\ell\in[1,n^N]:k_{\ell}>\sqrt{N}} |\alpha_{t,\ell}|=0$ if $k_{\ell}\leq\sqrt{N}$ for all $\ell\in[1,n^N]$.
The second sum on the right-hand-side of (\ref{bardtNdiffs*}) can clearly be upper-bounded by $n^ND/N$ by using (\ref{eq:s*appx3}).
The third sum on the right-hand-side of (\ref{bardtNdiffs*}) is upper-bounded by $\rho^N$ due to (\ref{eq:Dbound}).
As a result, we have, for each $t\in[0,T)$, $\delta>0$,
\begin{align*}
|\bar{d}_t^N - E^{\gamma^{*,N}} [\bar{d}_t^N] |   \leq &  \max_{\ell\in[1,n^N]:k_{\ell}>\sqrt{N}} |\alpha_{t,\ell}| +  n^N D/\sqrt{N} +  n^N D/N   + \rho^N.
\end{align*}
In view of (\ref{eq:alphabound}), this leads to, for each $t\in[0,T)$,  $\delta>0$,
\begin{align*}
& P^{s^{*,N}}\left[|\bar{d}_t^N - E^{\gamma^{*,N}} [\bar{d}_t^N] | > \delta + n^N D/\sqrt{N} +  n^N D/N   + \rho^N \right]  \\
& \quad \leq   P^{s^{*,N}}\left[\max_{\ell\in[1,n^N]:k_{\ell}>\sqrt{N}} |\alpha_{t,\ell}| > \delta \right]  \leq \frac{D^2}{\sqrt{N}\delta^2}.
\end{align*}
%\begin{align*}
%P^{s^{*,N}}[|\bar{d}_t^N - E^{\gamma^{*,N}} [\bar{d}_t^N ] | \geq 2\delta+n^N D (r+2/N)]   \leq &  1 - \bigg(1-\frac{D^2}{(rN-1) \delta^2}\bigg)^{n^N}.
%\end{align*}
Noting that $\lim_{N\rightarrow\infty} n^N D/\sqrt{N} +  n^N D/N   + \rho^N =0$, we obtain, for each  $t\in[0,T)$, $\delta>0$,
\begin{align}
\lim_{N\rightarrow\infty}P^{s^{*,N}}[|\bar{d}_t^N - E^{\gamma^{*,N}} [\bar{d}_t^N]| > 2 \delta]   = 0.
\label{eq:inequ1z}
\end{align}

Step (iii):
Regardless of how the joint strategy $s=(s^i)_{i\in[1,N]}$ in (\ref{eq:barstN}) is selected, we have $d_0^N=\bar{d}_0^N$. Therefore, by (\ref{eq:notcipdg}) and (\ref{eq:notcipds}), we have
\begin{equation}
|d_0^N - E^{\gamma^{*,N}}[\bar{d}_0^N] |\xrightarrow[N\rightarrow\infty]{\gamma^{*,N}}0  \qquad\textrm{and}\qquad |d_0^N - E^{\gamma^{*,N}}[\bar{d}_0^N] |\xrightarrow[N\rightarrow\infty]{s^{*,N}}0.
\label{eq:inequ2}
\end{equation}
By (\ref{eq:inequ2}) and the continuity of $g_0$, we have
\begin{align*}
m_{\mathbb{Y}} (y_0^N , \bar{y}_0^N ) \xrightarrow[N\rightarrow\infty]{\gamma^{*,N}} 0 \qquad\textrm{and}\qquad m_{\mathbb{Y}} (y_0^N , \bar{y}_0^N ) \xrightarrow[N\rightarrow\infty]{s^{*,N}} 0.
\end{align*}
 From the uniform continuity of $f_0^1$ (recall that $f_0^1=f_0^2=\cdots$) and the uniform equi-continuity of the mappings in $\mathbb{S}^1$, we obtain
 \begin{equation}
 \sup_{i\in[1,N]}\left|\begin{bmatrix} x_1^i \\ u_1^i\end{bmatrix}-\begin{bmatrix} \bar{x}_1^i \\ \bar{u}_1^i\end{bmatrix}\right| \xrightarrow[N\rightarrow\infty]{\gamma^{*,N}} 0
 \qquad\textrm{and}\qquad
 \sup_{i\in[1,N]}\left|\begin{bmatrix} x_1^i \\ u_1^i\end{bmatrix}-\begin{bmatrix} \bar{x}_1^i \\ \bar{u}_1^i\end{bmatrix}\right| \xrightarrow[N\rightarrow\infty]{s^{*,N}} 0.
 \label{eq:inequ2b}
 \end{equation}
 From (\ref{eq:inequ2b}) and  (\ref{eq:negdm1}), we also obtain
\begin{align*}
|d_1^N - \bar{d}_1^N | \xrightarrow[N\rightarrow\infty]{\gamma^{*,N}} 0  \qquad\textrm{and}\qquad |d_1^N - \bar{d}_1^N | \xrightarrow[N\rightarrow\infty]{s^{*,N}} 0
\end{align*}
which, due to (\ref{eq:notcipdg})-(\ref{eq:notcipds}), implies
\begin{equation}
|d_1^N - E^{\gamma^{*,N}} [\bar{d}_1^N] | \xrightarrow[N\rightarrow\infty]{\gamma^{*,N}} 0 \qquad\textrm{and}\qquad |d_1^N - E^{\gamma^{*,N}} [\bar{d}_1^N] | \xrightarrow[N\rightarrow\infty]{s^{*,N}} 0 .
\label{eq:inequ3}
\end{equation}
Now, from (\ref{eq:inequ2}) and (\ref{eq:inequ3}) and the continuity of $f_0^0$, $g_1$, we obtain
\begin{align*}
m_{\mathbb{Y}} (y_1^N,\bar{y}_1^N)  \xrightarrow[N\rightarrow\infty]{\gamma^{*,N}} 0 \qquad\textrm{and}\qquad m_{\mathbb{Y}} (y_1^N,\bar{y}_1^N)  \xrightarrow[N\rightarrow\infty]{s^{*,N}} 0.
\end{align*}
Continuing along the same lines gives us (\ref{eq:notcip})-(\ref{eq:notcip2}).

Step (iv):
We recall that $\gamma^{*,N}\in\Gamma$ is an SEBEU and, due to (\ref{eq:hatslopt}), $s^{*,N}$ satisfies
\begin{equation}
\label{eq:sebeup}
E^{\zeta_{\gamma^{*,N}}} [ J^i(s^{i,N},\cdot) ] \leq E^{(\tilde{\gamma}^i,\zeta_{\gamma^{*,N}})} [ J^i ], \quad \forall i \in[1,N], \tilde{\gamma}^i\in\Gamma^i.
\end{equation}
We also recall that $J^i=J^1$, $\Gamma^i=\Gamma^1$, for all $i\in\mathbb{N}$, where $J^1$ is continuous and bounded and  $\Gamma^1$ is compact.
As a result, (\ref{eq:notcip})-(\ref{eq:notcip2}) and (\ref{eq:sebeup}) imply that, for any $\epsilon>0$, there exists $\hat{N}_{\epsilon}\in\mathbb{N}$ such that, if $N\geq \hat{N}_{\epsilon}$, then
\begin{align*}
J^i (s^{i,N},\bar{y}^N) & \leq E^{\tilde{\gamma}^i} [J^i (\cdot,\bar{y}^N)] + \epsilon/2,  \quad \forall i \in[1,N], \tilde{\gamma}^i\in\Gamma^i \\
E^{\zeta_{s^{*,N}}} [J^i (s^{i,N},\cdot)] & \leq E^{(\tilde{\gamma}^i,\zeta_{s^{*,N}})} [J^i]  + \epsilon, \qquad   \forall i \in[1,N], \tilde{\gamma}^i\in\Gamma^i
\end{align*}
where $\bar{y}^N:=(\bar{y}_t^N)_{t\in[0,T)}$.
This concludes the proof.

\end{document}